\documentstyle[12pt]{article}
\def\~{\tilde}

\def\~{\tilde}

\def\cz{{\bf C}}
\def\nz{{\rm I \! N}}
\def\gz{{\rm Z \!\! Z}}
\def\pr{{\bf P}}

\def\r{\rightarrow}

\def\s{\subset}
\def\seq{\subseteq}
\def\se{\setminus}

\def\lr{\longrightarrow}

\def\co{ {\cal O} }
\def\ol{ \bar }

\def\da{ \downarrow }
\def\ua{ \uparrow}

\def\lra{ -\!\!-\!\!\!\longrightarrow }

\def\bpk{ \ol{P}_{k} }
\def\bpkl{ \ol{P}_{k + l} }

\def\bpl{ \ol{ P}_{l} }

\def\ph{\varphi}

\catcode`\"=\active
\def\3{\char\ss}

\def\qed{\ifmmode\sq\else{\unskip\nobreak\hfil
\penalty50\hskip1em\null\nobreak\hfil\sq
\parfillskip=0pt\finalhyphendemerits=0\endgraf}\fi}
\def\sq{\hbox{\rlap{$\sqcap$}$\sqcup$}}

\newtheorem{defi}{Definition}[section]
\newtheorem{prop}[defi]{Proposition}
\newtheorem{theo}[defi]{Theorem}
\newtheorem{lem}[defi]{Lemma}

\newtheorem{cor}[defi]{Corollary}

\newtheorem{mlem}[defi]{Main Lemma}

\begin{document}

\begin{center}
{\LARGE Logarithmic  Jet Bundles
and Applications}\\

\vspace{.7cm}
{\large Gerd-Eberhard Dethloff} and {\large Steven Shin-Yi Lu} \vspace{.3cm}
\end{center}

\setcounter{section}{-1}

\section{Introduction}
Hyperbolic complex manifolds have been studied
extensively during the last 30 years (see, for example, \cite{Ko1},
\cite{La1}).
However, it is still an
important problem in hyperbolic geometry to understand the algebro-geometric
and the differential-geometric meaning of hyperbolicity.
The use of jet bundles has become a powerful tool to attack this problem.
For example, Green and Griffiths (\cite{GG}) explained an approach to establish
Bloch's Theorem on the algebraic degeneracy of holomorphic maps
into abelian varieties by
constructing negatively curved pseudometrics on jet bundles and by
applying Ahlfors' Lemma. Siu and Yeung
(\cite{SY2}) clarified this approach. Moreover, they gave a Second Main Theorem
for divisors in abelian varieties, which was, very recently, clarified and
generalized to the case of semi-tori by Noguchi, Winkelmann and Yamanoi
(\cite{NWY}).

 Demailly  (\cite{Dem}) presented a new
construction of projective jets and pseudo-metrics on them
which realizes directly the approach to Bloch's theorem given in \cite{GG}.
 These projective jets are closer to the geometry of holomorphic curves
than the usual
jets, since the action of the group of reparametrizations of germs of curves,
which is geometrically redundant, is divided out. Using these pseudometrics on
projective jets,
Demailly and El Goul  (\cite{DemEG1}, see also Mc Quillan (\cite{MQ})) were
able to show that a (very) generic
surface $X$ in $\pr^3$ of degree $d \geq 21$ is Kobayashi hyperbolic. As a
corollary one obtains that
the complement of a (very) generic curve in $\pr^2$ of degree $d
\geq 21$
is hyperbolic and
hyperbolically embedded, a result  first proved
 by Siu and Yeung (\cite{SY1}) for much higher degree, using jet bundles and
value distribution theory. In both papers this quasi-projective case is
treated by
proving hyperbolicity of a branched cover over the compactification.

 However,
it is desirable to have also a direct approach to deal with quasiprojective
varieties, since one can hope to get easier proofs and even better results
\footnote{Actually, El Goul told the first named author that, using the results
of this paper, he succeeded to drop the degree in \cite{DemEG1} from 21 to
15.}.
So one should also consider the case of logarithmic jet bundles. Noguchi
(\cite{No1}) did  this  already for the case of the jet bundles used by
Green-Griffiths. Via these bundles he generalized Bloch's theorem to
semi-abelian
varieties.

The main purpose of the present
paper is to  generalize
Demailly's construction of projective jet bundles and strictly negatively
curved pseudometrics on them to the logarithmic case.
In sections 1 to 3, we establish this logarithmic generalization of
Demailly's construction explicitly via coordinates, just as Noguchi's
generalization of the jets used by Green-Griffiths. These explicit coordinates
 should be very useful for further applications.
We also have another, more intrinsic way to obtain the same generalization
in \cite{DL},
which is much shorter, but does not give coordinates right away.
In section 4
we  prove the Ahlfors Lemma and the Big Picard Theorem for logarithmic
projective jet bundles.

In section 5, we use our method to give a metric proof of Lang's Conjecture
for semi-abelian varieties and of a Big Picard analogue of it. The first result
is due to Siu and Yeung (\cite{SY}) and
Noguchi (\cite{PNog}), who used value distribution theory while
we use negatively curved jet metrics.
However, a common ingredient, due to Siu and Yeung (\cite{SY}), is
to construct a special jet differential (which naturally lives on the
jet space constructed by Demailly) from theta functions
on an abelian variety, the existence of which on a
semi-abelian variety we cite from
Noguchi (\cite{PNog}).
Hence, the main importance of this section is the method of proof. In
fact, we have to overcome some small technical difficulties to make our
method work in this case: For example, we have to introduce a d-operator
for sections over logarithmic projective jet bundles, and we have to deal
with the case of a divisor which can have worse singularities than normal
crossing, and with the precise relations between two different logarithmic
structures (the one coming from the boundary divisor of a semi-abelian
variety, the other coming from its union with an arbitrary reduced algebraic
divisor). In this way, section 5 can also serve as a complement
to sections 1 to 4.

We would like to thank J.P. Demailly and J. Noguchi for many
discussions on this subject.
We also would like to thank
the JSPS, the SFB, the DFG, the MSRI
 and the universities of G\"ottingen, of Waterloo and of Osaka
for their support during the preparation of this work. Finally we would
like to thank the referee for
his proposals, which led to a substantial improvement of this article.

\section{Log-directed jet bundles}
\subsection{Logarithmic jet bundles}

In this subsection we recall some basic
setup and results of  Noguchi in \cite{No1}. For the proofs
we refer to this article. Furthermore, in sections 1 to 3 we denote open
subsets
of a manifold by $O$, in order to distinguish them from open neighborhoods
of a given
point, usually with fixed coordinates centered at this point, which we
denote by $U$.\\

Let $X$ be a complex manifold.
Let $x \in X$. We consider germs $f:(\cz,0) \rightarrow  (X,x)$ of holomorphic
curves through $x$.  Two such germs are considered to be equivalent if they
have the same Taylor expansions of order $k$
in some (and hence, any) local coordinate around $x$. Denote the equivalence
class of $f$ by $j_k(f)$. We define
 $J_k  X_x = \{j_k(f) | f:(\cz,0) \rightarrow (X,x) \}$
and $J_k X = \bigcup_{x \in X} J_k  X_x$. Let $\pi :J_kX \r X$
be the natural projection. Then $J_kX$ carries the structure
of a holomorphic fiber bundle over $X$. It is called the {\it $k$-jet
bundle} over
$X$. If no confusion arises, we will denote the sheaf of sections of $J_kX$
also by
$J_kX$.
There exist, for $k \geq l$, canonical projection maps
\begin{equation} \label{tt0}
\pi_{l,k}:J_kX \rightarrow J_lX\; ; \;\; j_k(f) \r j_l(f),
\end{equation}
and $J_1X$ is canonically isomorphic to the holomorphic tangent bundle $TX$
over $X$.
 If $F: X \rightarrow Y$ is a holomorphic map
to another complex manifold $Y$, then it induces a holomorphic map
\begin{equation} \label{t0}
F_k:J_kX \rightarrow J_kY\; ; \;\; j_k(f) \r j_k(F \circ f)
\end{equation}
over $F$.

Let $\Omega X$ be the holomorphic cotangent bundle over $X$.
Take a holomorphic section $\omega \in H^0(O, \Omega X)$ for some open subset
$O \s X$. For $j_k(f) \in J_kX|_O$ we put $f^*\omega = Z(t)dt$. Then the
derivatives $\frac{d^j Z}{dt^j}(0)$, $0 \leq j \leq k-1$ are well defined,
independantly
of the representative $f$ for $j_k(f)$. Hence, we have a well defined mapping
\begin{equation} \label{t1}
\tilde{\omega}: J_kX|_O \r \cz^k\; ;
 \;\; j_k(f) \r (\frac{d^j Z}{dt^j}(0))_{0 \leq j \leq k-1}
\end{equation}
which is holomorphic.
If, moreover, $\omega^1$,...,$\omega^n$ with $n= \dim X$ are holomorphic
1-forms on $O$ such that $\omega^1 \wedge ... \wedge \omega^n$ does not
vanish anywhere, then we have a biholomorphic map
\begin{equation} \label{t2}
 (\tilde{\omega}^1,..., \tilde{\omega}^n) \times \pi: J_kX|_O \r
 (\cz^k)^n  \times O
\end{equation}
which we call the {\it trivialization} associated with
$\omega^1$,...,$\omega^n$.
More generally, if $\omega$ is a section over $O$ in the sheaf of
meromorphic 1-forms,
then the map $\tilde{\omega}$ defined as in equation (\ref{t1}) induces
a meromorphic vector valued function
\begin{equation} \label{t3}
\tilde{\omega}: J_kX|_O \r \cz^k.
\end{equation}

Let $\bar{X}$ be a complex manifold with a normal
crossing divisor $D$. This means that around any point $x$ of $\bar{X}$,
there exist local coordinates $z_1,...,z_n$
centered at $x$ such that $D$ is defined by
$z_1z_2...z_l=0$ in some neighborhood of $x$ and for some $l\leq n$.
We note that $l$ depends on $x$, which is implicitly assumed.
The pair $(\bar{X}, D)$ will be called a {\it log-manifold}.
Let $X= \bar{X} \setminus D$.

Following Iitaka (\cite{Ii}), we define the logarithmic
cotangent sheaf $\bar{\Omega} X
= \Omega(\bar X, {\rm log} D)$ as the locally free
subsheaf of the sheaf
of meromorphic 1-forms on $\bar{X}$,
 whose restriction to $X$ is $\Omega X$
(where we identify vector bundles and their sheaves of sections)
and whose localization at any point $x\in D$ is given by
\begin{equation} \label{t4}
\bar{\Omega} X_x=
\sum_{i=1}^l {\cal O}_{\bar X,x}{dz_i\over z_i} +
\sum_{j=l+1}^n {\cal O}_{\bar X,x} dz_j,
\end{equation}
where the local coordinates $z_1,...,z_n$ around $x$ is chosen as before.
Its dual, the logarithmic tangent sheaf $\bar{T}X =T(\bar X, - {\rm log} D)$,
is a locally free subsheaf of
the holomorphic tangent bundle $T \bar{X}$ over $\bar{X}$.
Its restriction to $X$ is identical to $T X$,
and its localization at any $x\in D$ is given by
\begin{equation} \label{t5}
\bar{T} X_x=
\sum_{i=1}^l {\cal O}_{\bar X,x}z_i{\partial \over \partial z_i} +
\sum_{j=l+1}^n {\cal O}_{\bar X,x} {\partial \over \partial z_j}.
\end{equation}
 Given
log-manifolds  $(\bar{X'},D')$ and $(\bar{X},D)$,
a holomorphic map $F:\bar{X'} \r \bar{X}$ such that $F^{-1}D \s D'$ will be
called a {\it log-morphism} from $(\bar{X'},D')$ to
$(\bar{X},D)$. If no confusion arises, we will simply write
 $F: X' \r X$ for the log-morphism
$F:(\bar{X'},D') \r (\bar{X},D)$.
It induces (see \cite{Ii}) vector bundle morphisms,
\begin{equation} \label{tt3}
F^*: \bar{\Omega} X \r F^{-1}\bar{\Omega} X' \r \bar{\Omega} X'\ \mbox{and}\
F_*: \bar{T} X' \r F^{-1}\bar{T} X \r \bar{T} X ,
\end{equation}
where we have again
identified locally
free sheaves and vector bundles.\\

Let $s \in H^0(O, J_k \bar X)$ be a holomorphic section over an open
subset $O \s \bar X$. We say that $s$ is a {\it logarithmic $k$-jet field}
if the map $\tilde{\omega} \circ s|_{O'} : O' \r \cz^k$ is holomorphic
for all $\omega \in H^0(O', \bar \Omega X)$ for all open subsets $O'$ of
$O$ and where the map $\tilde{\omega}$ is defined as in equation
(\ref{t3}).
The set of logarithmic $k$-jet fields over open subsets of $\bar X$ defines a
subsheaf
of the sheaf $J_k \bar X$, which we denote by $\bar J_kX$. By a) of the
following
proposition, $\bar J_k X$ is the sheaf of
sections of a holomorphic fiber bundle
over $\bar X$, which we denote again by $\bar J_kX$, and which we call the
{\it logarithmic $k$-jet bundle} of $(\bar X,D)$.
\begin{prop}[see \cite{No1}]
\label{p1} $ $\\[-.3in]
\begin{itemize}
\item[a)] $\bar J_k X$ is the sheaf of
sections of a holomorphic fiber bundle
over $\bar X$.
   (However, it is only a subsheaf and not a subbundle of $J_k \bar X$.)
\item[b)]   We have a canonical identification of $(\bar J_k X)|_X$ with
$J_k X$.
\item[c)] Let $O \s  \bar X$ be an open set and $\theta$ be any meromorphic
function
   on $O$ such that the support of its divisor $(\theta )$  is contained in
$D$.
   Let $d^l \log \theta$ be the
   $l$-th
   component of the map $\widetilde{\Theta}: \bar J_k X|_O \r \cz^k$,
   where $\Theta = d \log \theta$
   (see equation (\ref{t1}) and (\ref{t3})).
   Then the differentials $d^l \log \theta$, $l=1,...,k$, define holomorphic
   functions on $\bar J_k X|_O$.
   Moreover, outside the support of $(\theta )$, we have
   $(d^l \log \theta )(j_k(f)) =
   \frac{\partial ^l \log (\theta \circ f)}{\partial t^l}(0)$.
\item[d)] There exists, for $k \geq l$, a canonical projection map
   $\pi_{l,k}:\bar J_kX \rightarrow \bar J_lX$,
   which extends the map
   $(\pi_{l,k}|_{J_kX}): J_kX \r J_lX$ (see equation (\ref{tt0})),
   and $\bar J_1X$ is canonically isomorphic to $\bar TX$.
\item[e)] A log-morphism  $F: X' \r X$ induces a canonical map
   $F_k:\bar J_kX' \rightarrow \bar J_kX$, which extends the map
   $F_k|_{J_kX}: J_kX' \r J_kX$ (see equation (\ref{t0})).
\end{itemize}
\end{prop}

Finally, we want to express the local triviality of $\bar J_k X$ explicitly
in terms of coordinates. Let $z_1,...,z_n$ be coordinates in an open set $U
\s \bar X$
in which $D=\{z_1z_2...z_l=0\}$.
Let $\omega^1= {dz_1\over z_1},...,\omega^l= {dz_l\over z_l},\;
\omega^{l+1}= dz_{l+1},...,\omega^{n}= dz_{n}$.
Then we have a biholomorphic map (see equations (\ref{t2})
and (\ref{t3}))
\begin{equation} \label{t6}
  (\tilde{\omega}^1,..., \tilde{\omega}^n) \times \pi: \bar J_kX|_U \r
  (\cz^k)^n \times U.
\end{equation}
Let $s \in H^0(U, \bar J_kX)$ be given by
$s(x)=(Z(x);x)$  in this trivialization with
$$Z =(Z^i_j)_{i=1,...,n; j=1,...,k} , $$
where the $Z_j^i(x)$ are holomorphic functions on $U$ and the indices $j$
correspond to the orders of derivatives. Then the same $s$, considered
as an element of $ H^0(U, J_k \bar X)$ and trivialized by
$\omega^1=dz_1,...,\omega^n=dz_n$ (see equation (\ref{t2})) is given by
$s(x)=(\hat{Z}(x);x)$ with
$\hat{Z} =(\hat{Z}^i_j)_{i=1,...,n; j=1,...,k}$,
where
\begin{equation} \label{t7}
\hat{Z}^i_j=\left\{
\begin{array}{ccc}
z_i(Z^i_j + g_j(Z^i_1,...,Z^i_{j-1})) & : & i \leq l\\
Z^i_j & : & i \geq l+1
\end{array} \right. .
\end{equation}
Here, the $g_j$ are polynomials in the variables $Z^i_1,...,Z^i_{j-1}$ with
constant coefficients  and without constant
terms (in particular $g_1=0$), which are
obtained by expressing first the different components $Z_j^i$ of
$\widetilde{(\frac{dz_i}{z_i})} \circ s(x)$ in terms of the components
$\hat{Z}^i_j$
of $\widetilde{dz_i} \circ s(x)$ by using the chain rule, and then by
expressing the
$Z_j^i$ in terms of the $\hat{Z}^i_j$ by
inverting this system of polynomial equations.
This clarifies equation (1.13) in \cite{No1}, where, for $i \leq l$, only
the leading term $z_iZ^i_j$ is given.
This also exhibits the sheaf inclusion
$\bar J_kX|_U \s J_k \bar X|_U$ explicitly in terms of coordinates.
By abuse of notation, we also consider the $Z^i_j$'s as the holomorphic
functions defined on $\bar J_kX|_U$ given by equation (\ref{t6}), so that
$Z_1^1,...,Z^n_k;z_1,...,z_n$ form a holomorphic coordinate system on $\bar
J_kX|_U$.

We remark that a trivialization of $\bar J_k X|_U$ is also obtained
if we replace the special $\omega$'s used in equation (\ref{t6})
by any $\omega^1,...,\omega^n \in H^0(U, \bar \Omega X)$ with
\begin{equation} \label{t8}
\omega^1 \wedge ... \wedge \omega^n =\frac{a(x)}{z_1z_2...z_l}
dz_1 \wedge ... \wedge dz_n,
\end{equation}
where $a(x)$ is a nowhere vanishing holomorphic function on $U$.

\subsection{Log-directed jet bundles}
We first follow Demailly (\cite{Dem}). Let $X$ be
a complex manifold together with a holomorphic subbundle
$V \s TX$. The pair
$(X,V)$  is called a {\it  directed
manifold}. If $(X,V)$ and $(Y,W)$ are two
such manifolds,
then a holomorphic map $F:X \r Y$ which satisfies $F_*(V) \s W$ is called a
{\it directed morphism}.

Let  $(X,V)$ be a directed manifold. The subset $J_kV$ of $J_kX$  is
defined to be the set of $k$-jets $j_k(f) \in J_kX$
for which there exists a representative
$f:(\cz,0) \r (X,x)$ such that $f'(t) \in V_{f(t)}$ for all $t$ in a
neighborhood of $0$. We will show in the next subsection that
$J_kV$ is a fiber bundle over $X$,
which we call the {\it directed k-jet bundle}
$J_kV$ of  $(X,V)$.
If $F:(X,V) \r (Y,W)$ is a directed morphism,
then equation (\ref{t0}) induces a holomorphic map
\begin{equation} \label{t9}
F_k:J_kV \rightarrow J_kW\; ; \;\; j_k(f) \r j_k(F \circ f)
\end{equation}
over $F$, since the restriction of $F_k :J_kX \r J_kY$ to $J_kV$
maps to $J_kW$ as $(F \circ f)'(t)=F_*(f'(t)) \in W_{F \circ f(t)}$ if
$f'(t) \in V_{f(t)}$.\\

We now generalize Demailly's directed $k$-jet bundles to the
logarithmic
context. We define a
{\it log-directed manifold}
to be the triple $(\bar X,D,\bar V)$, where $(\bar X,D)$ is a log-manifold
together with a subbundle $\bar V$ of
$\bar T X$. A {\it log-directed morphism}
between log-directed manifolds
$(\bar{X'},D',\bar V')$ and $(\bar{X},D,\bar V)$ is
a log-morphism $F :(\bar X',D') \r (\bar X, D)$
such that $F_*(\bar V') \s \bar V$.

Let $(\bar X,D,\bar V)$ be a log-directed manifold and set $V=\bar V|_X$.
By Proposition \ref{p1} we can canonically identify $(\bar J_k X)|_X$ with
$J_kX$. Hence,
the directed $k$-jet
bundle $J_k V$ of $(X, V)$  can be considered as a subset of the logarithmic
$k$-jet bundle $\bar J_k X$ over $\bar X$. We define the {\it log-directed
$k$-jet bundle} $\bar J_kV$ of $(\bar X, D, \bar V)$ to be the topological
closure
$\overline{J_k V} \s \bar J_k X$ of $J_kV$ in $\bar J_kX$.
If $F:(\bar X', D',\bar V') \r (\bar X, D, \bar V)$ is a log-directed
morphism,
it induces a  map
\begin{equation} \label{t10}
F_k:\bar J_kV' \rightarrow \bar J_kV
\end{equation}
over $F$ which is holomorphic. It is the restriction of the canonical map
 $F_k:\bar J_kX \rightarrow \bar J_kX'$ (see Proposition \ref{p1}) to $\bar
J_kV'$ and is also an
 extension of the map
$F_k|_X: J_kV' \rightarrow  J_kV$ (see equation (\ref{t9})) to $\bar J_kV'$.

\subsection{Structure of log-directed jet bundles}
In this subsection, we study the local structure of
$\bar J_kV \s  \bar J_kX$ over $\bar X$. In particular, we show
that $\bar J_kV \s  \bar J_kX$ is a submanifold of  $\bar J_kX$ which itself is
a locally trivial bundle.
This justifies the name of log-directed $k$-jet bundle for $\bar J_kV$
introduced in the
previous subsection.\\

First, we consider the directed manifold $(X,V)$. For any point $x_0 \in X$,
 there is a
coordinate system $(z_1,...,z_n)$, centered at $x_0$, on a neighborhood $U$ of
$x_0$ such that the fibers
$V_x$ for $x \in U$
can be defined by linear equations
\begin{equation} \label{t11}
V_x = \{ \xi = \sum_{1 \leq i \leq n} \xi_i \frac{\partial}{\partial z_i} \; |
\;\; \xi_i = \sum_{1 \leq m \leq r} a_{im}(x) \xi_m \; {\rm for} \;
i=r+1,...,n \}\, .
\end{equation}
We fix $x_0$, $U$ and these coordinates from now on.
If we trivialize $TX=J_1X$ by $\omega^1= dz_1,...,\omega^n=dz_n$ over $U$
as in equation (\ref{t2}),
we obtain
\begin{equation} \label{t11a}
V_x = \{(Z_1^1,...,Z_1^n) \; | \;\;
 Z_1^i = \sum_{1 \leq m \leq r} a_{im}(x) Z_1^m \; {\rm for} \;
 i=r+1,...,n \} \, .
\end{equation}
If we trivialize $J_kX$ by the same forms,
we obtain more generally:
\begin{prop} \label{J_kV} $ $\\[-.3in]
\begin{itemize}
\item[a)] Let $P_h^i$ be the polynomials in the variables $Z_j^i$
with coefficients depending holomorphically on $x$
obtained by formally differentiating the  equations
$$f_i'(t) = \sum_{1 \leq m \leq r} a_{im}(f(t)) f_m'(t)$$
$h-1$ times with respect to $t$,  using the fact that
$Z^i_j(j_k(f))= f_i^{(j)}(t)$ in our trivialization. Then we have
\begin{equation} \label{t12}
(J_kV)_x =\{ (Z^i_j)_{i=1,...,n; j=1,...,k}
 \; | \; \;
 Z_h^i = P_h^i(x,Z_1^1  , ... , Z_1^n,...,
Z_{h-1}^1  , ... $$
$$..., Z_{h-1}^n,Z_h^1  , ... , Z_h^r)
  \; {\rm for} \;
 h=1,...,k, \; \; i=r+1,...,n \}.
\end{equation}

\item[b)]
 $J_kV \s J_kX$ is a submanifold, and the canonical
projection
$$K: U \r \cz^r\; ; \;\; (z_1,...,z_n) \r (z_1,...,z_r) $$
induces a bundle isomorphism
$$K_k: J_kV|_U \r K^{-1}(J_k\cz^r).$$
\end{itemize}
\end{prop}
\noindent
{\bf Proof for a) } Let $j_k(f) \in J_kV$. By definition,
 there exists a
representative
$f$ such that $f'(t) \in V_{f(t)}$
for all $t$ in a neighborhood of $0 \in \cz$, namely
$$f_i'(t) = \sum_{1 \leq m \leq r} a_{im}(f(t)) f_m'(t)\, .$$
Now it follows from the chain rule that $j_k(f)$ satisfies equations of the
form
$Z_h^i = P_h^i$, $h=1,...,k, \;  i=r+1,...,n$, given in equation (\ref{t12}).

Conversely let $Z^i_j \in \cz$, $i=1,...,n, \; j=1,...,k $ be given
satisfying the equations
$Z_h^i = P_h^i$, $h=1,...,k, \;  i=r+1,...,n$ of equation (\ref{t12}).
For $x \in X$ fixed, define, for $i=1,...,r$, holomorphic functions
$$f_i: \cz \r \cz; \; t \r z_i(x) + \sum_{\nu =1}^k
\frac{Z^i_{\nu}}{\nu !} t^{\nu}.$$
Now we integrate the system of differential equations
$$f_i'(t) = \sum_{1 \leq m \leq r} a_{im}(f_1(t),...,f_n(t)) f_m'(t)\;\;
 i=r+1,...,n\, ,$$
to obtain a germ $f:(\cz,0) \r (X,x)$ with $z_i(f)=f_i$, $i=1,...,n$.
We see by construction (as $f'(t) \in V_{f(t)}$) that
$$(\tilde{\omega}^1,..., \tilde{\omega}^n)(j_k(f))=
(Z^i_j)_{i=1,...,n; j=1,...,k} \:.$$
{\bf Proof for b) } If one replaces successively in the $P^i_h$ all the
$Z^i_j$ with $i \geq r+1$ and $j \leq h-1$ by their expressions in terms
of the $Z^i_j$ with $i \leq r$ via equation (\ref{t12}), we get from  a)
that $J_kV|_U$ is the graph
of these new functions
$P_h^i$, $i=r+1,...,n; \; h=1,...,k$ in the variables $z_1,...,z_n$ and
$Z^i_j$, $i=1,...,r; \; h=1,...,k$. These are in turn coordinates for
$K^{-1}(J_k \cz ^r)$.  \qed
\medskip

Let now $(\bar X, D, \bar V)$ be a log-directed manifold.
Let $x_0 \in \bar X$ and let $z_1,...,z_n$ be a coordinate system
centered at $x_0$ on a neighborhood $U$ of $x_0$ where $D$
is defined by $z_1z_2...z_l=0$ for some $l \leq n$.
If we trivialize $\bar TX = \bar J_1X$ over $U$ by
 $\omega^1= {dz_1\over z_1},...,\omega^l= {dz_l\over z_l},\;
\omega^{l+1}= dz_{l+1},...,\omega^{n}= dz_{n}$
as in equations (\ref{t6}) and (\ref{t7}), we obtain
\begin{equation} \label{t13}
V_x = \{(Z_1^1,...,Z_1^n) \; | \;\;
 Z_1^i = \sum_{m \in A} a_{im}(x) Z_1^m \; {\rm for} \;
 i \in B  \}
\end{equation}
for all $x\in U$,
where, after permuting $z_1,...,z_l$ respectively $z_{l+1},...,z_n$,
we have $A= \{ 1,...,a,l+1,...,l+b \}$ and $B=\{ 1,...,n \} \setminus A$
with $a+b=r= {\rm rank}\, V$. We fix this setup for the rest of this section.

First, we generalize the projection $K$ to log-directed manifolds:
\begin{prop} \label{pro}
With $E= \{ z_1...z_a=0 \}$, the log-directed projection
$$K: (\bar X,D,\bar V)|_U \r  (\cz^r, E, \bar T \cz^r) \; ; \; \;
(z_1,..., z_n) \r (z_1,...z_a,z_{l+1},...,z_{l+b})$$
has bijective differential map $(K_*)_x$ for all  $x \in U$.
 \end{prop}
\noindent {\bf Proof }
We trivialize $\bar T \cz^r$ by the forms
$\omega^1= {dz_1\over z_1},...,\omega^a= {dz_a\over z_a},\;
\omega^{l+1}= dz_{l+1},...,\omega^{l+b} = dz_{l+b}$.
We claim that $K_*$ is given by the projection map
\begin{equation} \label{pf1}
(K_*)_x: (\bar T X)_x \r (\bar T \cz^r)_x; \;
(Z_1^1,...,Z_1^n) \r (Z_1^1,...,Z_1^a,Z_1^l,..., Z_1^{l+b})\,
\end{equation}
in these coordinates.
In fact, by analytic continuation it suffices to prove equation (\ref{pf1})
for $x \in X = \bar X \se D$. Let $(Z_1^1,...,Z_1^n) \in (\bar T X)_x = (T
X)_x$
be a vector in the logarithmic coordinate system.
If we retrivialize $(T X)_x$ respectively $(T\cz^r)_{K(x)}$ with the forms
$dz_i$ ($i=1,...,n$
respectively $i \in A$) instead, then the given vector is expressed by
 $(z_1Z_1^1,...,z_lZ_1^l, Z_1^{l+1},...,Z_1^n)$ (see equations (\ref{t6})
and (\ref{t7})). Furthermore, in the latter trivialization,
the map $(K_*)_x$ is just the projection to the components given by $A$.
So equation (\ref{pf1}) follows.
Hence, the assertion follows from equation (\ref{t13}).
\qed \medskip

If we trivialize $\bar J_kX$ by the same forms as in equation (\ref{t13}),
we obtain the following extension of Proposition \ref{J_kV}.
\begin{prop} \label{barJ_kV}
Let the setup be as above. \\[-.3in]
\begin{itemize}
\item[a)]
There are polynomials $Q_h^i$ in the variables $Z_j^i$
with coefficients which are holomorphic functions on $U$ such that
\begin{equation} \label{t14}
\bar J_kV_x =\{
(Z^i_j)_{i=1,...,n; j=1,...,k}
 \; | \; \;
 Z_h^i = Q_h^i(x,Z_1^1  , ... , Z_1^n,...,
Z_{h-1}^1  , ... $$
$$..., Z_{h-1}^n,Z_h^1  , ... , Z_h^r)
  \; {\rm for} \;
 h=1,...,k, \; \; i \in B \}.
\end{equation}

\item[b)]
$\bar J_kV \s \bar J_kX$ is a submanifold and the
projection map $K$ defined in Proposition \ref{pro}
induces a bundle isomorphism
$$K_k: \bar J_kV|_U \r K^{-1}(\bar J_k (\cz^r \se E)).$$
\end{itemize}
\end{prop}
\noindent {\bf Proof for a) }
By using the coordinate system $(z_1,...,z_n)$ on $U$, the map
$$\Psi : \cz^n \r \cz^n; \; (w_1,...,w_n) \r (e^{w_1},...,e^{w_l},
w_{l+1},...,w_n)
=(z_1,...,z_n)$$
induces a locally
biholomorphic
map $\Psi: \Psi^{-1} (U) \r U \se D$. Let
$\hat{V}=\Psi^{-1}_*(V) \s T(\Psi^{-1}(U))$
and let $W_j^i$, $i=1,...,n; \; j=1,...,k$ be the components  of
the first part of the trivialization map
$$ (\tilde{dw_1},..., \tilde{dw_n}) \times  \pi: J_k (\Psi^{-1} (U)) \r
  (\cz^k)^n \times \Psi^{-1} (U).$$
\begin{lem}\label{lemW}
On $J_k (\Psi^{-1} (U))$, we have
$$W_j^i = Z_j^i \circ \Psi_k \;\; i=1,...,n; \; j=1,...,k.$$
\end{lem}
{\bf Proof of Lemma \ref{lemW} } Let $j_k(f) \in J_k (\Psi^{-1}(U))$ and let
$f=(f_1,...,f_n): (\cz, 0) \r \Psi^{-1} (U)$ represent it.
We put $f^* dw_i = df_i(t) = C_i(t)dt$.
Then we have $W_j^i (j_k(f)) =
 \frac{\partial^{j-1}C_i(t)}{\partial t^{j-1}}|_{t=0}$.
On the other hand,
$$(\Psi \circ f)^* \omega^i = df_i(t) = C_i(t)dt$$
independently of $i$,
 and hence,
$$Z_j^i \circ \Psi_k (j_k(f)) =
\frac{\partial^{j-1}C_i(t)}{\partial t^{j-1}}|_{t=0}
= W_j^i (j_k(f)).$$ \\[-.5in] \qed\medskip

Since $\hat{V} \s T(\Psi^{-1} (U))$ is the inverse image of $V \s T(U \se D)$,
 we have
\begin{equation} \label{pf2}
\hat{V}_w = \{(W_1^1,...,W_1^n) \; | \;\;
 W_1^i = \sum_{m \in A} a_{im} \circ \Psi (w) \cdot W_1^m \; {\rm for} \;
 i \in B  \},
\end{equation}
for $w \in \Psi^{-1} (U)$. Using Proposition \ref{J_kV}, we get that
$$(J_k\hat{V})_w =\{
(W^i_j)_{i=1,...,n; j=1,...,k}
 \; | \; \;
 W_h^i = P_h^i(w,W_1^1  , ... , W_1^n,...,
W_{h-1}^1  , ... $$
$$..., W_{h-1}^n,W_h^1  , ... , W_h^r)
  \; {\rm for} \;
 h=1,...,k, \; \; i=r+1,...,n \},
$$
where the $P_h^i$ are the polynomials in the variables $W_j^i$
with coefficients depending holomorphically on $w \in \Psi^{-1}(U)$
obtained by formally differentiating
\begin{equation} \label{pf10}
f_i'(t) = \sum_{1 \leq m \leq r} a_{im}\circ \Psi(f(t)) \cdot f_m'(t)
\end{equation}
$h-1$ times. The important point is now that the coefficient functions factor
through $\Psi$ by holomorphic functions which are still holomorphic for all
$x \in U$:
\begin{mlem}
The coefficients of the polynomials $P_h^i$
factor through $\Psi$ by holomorphic functions
which are defined on all of $U$. Namely,
$$P_h^i(w,W_1^1  , ... , W_1^n,...,
W_{h-1}^1  , ..., W_{h-1}^n,W_h^1  , ... , W_h^r)$$ $$=
Q_h^i(x,W_1^1  , ... , W_1^n,...,
W_{h-1}^1  , ..., W_{h-1}^n,W_h^1  , ... , W_h^r)$$
for $x=\Psi (w)\in U\setminus D$,
where the $Q_h^i$ are  polynomials in the variables $W_j^i$,
with coefficients which are holomorphic in $x$ on  all of $U$.
\end{mlem}
{\bf Proof of the Main Lemma } If $\alpha : U \r \cz$ is holomorphic,
we have:
$$\frac{\partial}{\partial t} \alpha \circ \Psi (f(t)) =
\sum_{\mu = 1}^{n} (\frac{\partial \alpha}{\partial z_{\mu}} \circ \Psi)(f(t))
\cdot \frac{\partial}{\partial t} (\Psi_{\mu} \circ f)(t) $$
$$ = \sum_{\mu = 1}^{l} (\frac{\partial \alpha}{\partial z_{\mu}} \circ
\Psi)(f(t))
\cdot \frac{\partial}{\partial t} e^{f_{\mu}}(t) +
\sum_{\mu = l+1}^{n} (\frac{\partial \alpha}{\partial z_{\mu}} \circ
\Psi)(f(t))
\cdot \frac{\partial}{\partial t} f_{\mu}(t) $$
$$ = \sum_{\mu = 1}^{l} (\frac{\partial \alpha}{\partial z_{\mu}} \circ
\Psi)(f(t))
\cdot (z_{\mu} \circ \Psi) (f(t)) \cdot f_{\mu}'(t) +
\sum_{\mu = l+1}^{n} (\frac{\partial \alpha}{\partial z_{\mu}} \circ
\Psi)(f(t))
\cdot  f_{\mu}'(t) $$
$$ = \sum_{\mu = 1}^{l} ((z_{\mu} \frac{\partial \alpha}{\partial z_{\mu}})
\circ \Psi)(f(t))
\cdot  f_{\mu}'(t) +
\sum_{\mu = l+1}^{n} (\frac{\partial \alpha}{\partial z_{\mu}} \circ
\Psi)(f(t))
\cdot  f_{\mu}'(t) .$$
Now the assertion follows by induction on $h$ as the
coefficients of the
polynomials $P_h^i$ are obtained by formally differentiating the equation in
equation (\ref{pf10}) $h-1$ times and that the functions $a_{im}$ are
holomorphic on all of $U$. \qed

Finally we patch together these results and obtain the proof  of
Proposition \ref{barJ_kV} (a): Using  Lemma \ref{lemW}, the Main Lemma and
the local isomorphisms
$\Psi^{-1}$ we see that this assertion holds for all $x \in U \se D$,
with equations $Z_h^i = Q_h^i(\Psi ( \Psi^{-1}(x)),...W_j^i...) =
Q_h^i(x,...Z_j^i...)$ which are independent of the choice of
the local isomorphism $\Psi^{-1}$. Moreover, their coefficient functions are
still holomorphic on $U$. Since $\bar J_kV$ is defined as the closure of
$J_kV$ in $\bar J_k X$, the structure of the equations
$Z_h^i = Q_h ^i$ implies the assertion for all $x \in U$.\\
{\bf Proof for Proposition \ref{barJ_kV} b) }
It is verbatim that of Proposition \ref{J_kV} b).
\qed

\subsection{Regular jets}
Let $(X,V)$ be a directed manifold.
The subset $J_kV^{\rm sing} \s J_kV$ of {\it singular k-jets}
is defined to be the subset of $k$-jets $j_k(f) \in J_kV$ of germs
$f:(\cz,0) \r (X,x)$ such that $f'(0) =0$.
Its complement
$J_kV^{\rm reg} =J_kV \setminus J_kV^{\rm sing}$
defines the {\it regular k-jets}.

Let now $(\bar X,D, \bar V)$ be a log-directed manifold.
Define $\bar J_kV^{\rm sing} \s \bar J_kV$ to be the closure
$\overline{J_kV^{\rm sing}} \s  \bar J_kV$ of $J_kV^{\rm sing}$
in  $\bar J_kV$ and set
$\bar J_kV^{\rm reg} =\bar J_kV \setminus \bar J_kV^{\rm sing}$.

\begin{prop} \label{p} $ $\\[-.3in]
\begin{itemize}
\item[a)]
$\bar J_kV^{\rm sing} \s \bar J_kV$ is a smooth submanifold of codimension
$r= {\rm rank} \bar V$.
In terms of the local coordinates
of $\bar J_kV \s \bar J_kX$ (see Proposition \ref{barJ_kV}), this
submanifold is given by the equations
$$Z_1^i=0, i \in A.$$
\item[b)] The bundle isomorphism $K_k$ given in Proposition \ref{barJ_kV} b)
respects the singular and regular jets.
\end{itemize}
\end{prop}
\noindent {\bf Proof for a) } Using equations (\ref{t6}), (\ref{t7}) and
(\ref{t13}), we see that
$J_kV^{\rm sing}$ in  $\bar J_kV$ is given locally by the equations
$Z_1^i=0, i \in A$. So the assertion follows.\\
{\bf Proof for b) } This follows directly from the proof of
Proposition \ref{barJ_kV} b).
 \qed

\section{Log-Demailly-Semple jet bundles}
\setcounter{equation}{0}
A natural notion of higher order contact structures was
introduced on a firm setting by Demailly (\cite{Dem}) in
the holomorphic category for the study of complex hyperbolic
geometry. These structures realize natural
``quotient'' spaces of directed jet bundles.
Demailly called them Semple Jet bundles after Semple
(\cite{Sem}), who constructed and
worked with these bundles over ${\pr}^2$. In this section, we
generalize these bundles
to the logarithmic case and prove some important properties, like
functoriality and local triviality. Their connection with
log-directed jet bundles will be discussed in section 3.

\subsection{Definition of log-Demailly-Semple jet bundles}

We begin with a log-directed manifold $X_0=(\bar X_0, D_0, \bar V_0)$.
We inductively define $(\bar X_{k}, D_k, \bar V_k)$ as
follows. Let $\bar X_{k}={\pr}(\bar V_{k-1})$ with its natural projection
$\pi_{k}$ to $\bar X_{k-1}$. Set $D_{k}=\pi_k^{-1} (D_{k-1})$ and
$X_k=\bar X_k\se D_k$.
Let ${\cal O}_{\bar X_k}(-1)$ be the tautological
subbundle of $\pi_k^{-1}\bar V_{k-1}\seq \pi_k^{-1}\bar T X_{k-1}$, and set
\begin{equation}\label{O1}
\bar V_k=(\pi_k)_{\star}^{-1}\Big( {\cal O}_{\bar X_k}(-1)\Big).
\end{equation}
Equivalently, $\bar V_k \s \bar T X_k$ is defined,
for every point $(x, [v]) \in  \bar X_{k}={\pr}(\bar V_{k-1})$ associated with
a vector $v \in (\bar V_{k-1})_x$ for $x \in\bar X_{k-1}$, by
$$( \bar V_k)_{(x, [v])}= \{ \xi \in (\bar T X_k)_{(x,[v])}: \; (\pi_k)_* \xi
\in \cz v \}, \; \;\; \cz v \s ( \bar V_{k-1})_x \s (\bar TX_{k-1})_x.$$
Since $(\pi_k)_{\star}: \bar TX_k\r \pi_k^{-1}\bar T X_{k-1}$
has maximal rank everywhere as it is a bundle projection,
we see that $\bar V_k$
is a subbundle of $\bar TX_k$ giving a log-directed structure
for $X_k$ and also for $\pi_k$, thus completing our inductive definition.

We set $\bar P_k V=\bar X_k$, $ P_k V= X_k$,
 $\bar P_k X=\bar P_k TX$  and $ P_k X= P_k TX$.
Let
$$\pi_{j,k}=\pi_{j+1}\circ\cdots\circ\pi_{k-1}\circ\pi_{k}:
\bar P_k V\r \bar P_{j}V$$ for $j<k$. We also put
$(\bar P_kV)_x = (\pi_{0,k})^{-1}(x)$ and
$(\bar V_k)_x = \bar V_k|_{(\bar P_kV)_x}$ for $x \in \bar X$.

Note that
$\ker\, (\pi_k)_{\star}=T_{\bar P_kV/ \bar P_{k-1}V}$ by definition.
This gives the following short exact sequence
of vector bundles over $\bar P_kV$:
\begin{equation}
0\lr T_{\bar P_kV/ \bar P_{k-1}V}\lr \bar V_k
\stackrel{(\pi_k)_{\star}}{\lra} {\cal O}_{\bar P_kV}(-1)\lr 0.
\end{equation}
Furthermore, we have  the Euler
exact sequence for projectivized bundles (applied to the bundle
${\pr}(\bar V_{k-1}) \r \bar P_{k-1}V = \bar X_{k-1}$)
\begin{equation}
0\lr {\cal O}_{\bar P_kV}\lr \pi_k^{-1}\bar V_{k-1}\otimes {\cal O}_{\bar
P_kV}(1)
\lr T_{\bar P_kV/ \bar P_{k-1}V}\lr 0.
\end{equation}
 The composition of vector bundle morphisms over $\bar P_kV$
$${\cal O}_{\bar P_{k}V}(-1)\hookrightarrow \pi_k^{-1} \bar V_{k-1}
\stackrel{(\pi_k)^{-1}(\pi_{k-1})_{\star}}{\lra}
\pi_k^{-1}{\cal O}_{\bar P_{k-1}V}(-1)$$
yields an effective divisor $\Gamma_k$ corresponding to a section of
\begin{equation}
{\cal O}_{\bar P_{k}V}(1)\otimes \pi_k^{-1}{\cal O}_{\bar P_{k-1}V}(-1)
={\cal O}(\Gamma_k).\end{equation}
  There is a canonical divisor on $\bar P_kV$ given by
$$\bar P_kV^{\rm sing}=\bigcup_{2\leq j\leq k} \pi_{j,k}^{-1}(\Gamma_j)
\subset \bar P_k V.$$ Finally, we set
$$\bar P_kV^{\rm reg} = \bar P_k V\se \bar P_kV^{\rm sing}\ \ \mbox{ and }
\ \ {\cal O}_{ \bar P_k V}(-1)^{\rm reg}=
\Big({\cal O}_{ \bar P_k V}(-1)|_{\bar P_kV^{\rm reg} }\Big)\se
\bar P_k V,$$
where the last $\bar P_k V$ denotes the zero section.

\subsection{Properties of log-Demailly-Semple jet bundles}

\begin{prop} \label{functoriality}
Let $F\colon (\bar X', D',\bar V')\! \r \!(\bar X, D, \bar V)$
be a log-directed morphism.\\[-.25in]
\begin{itemize}
\item[a)] For all $k \geq 0$ there exist log-directed meromorphic maps
(log-directed morphisms outside the locus of indeterminacy)
$$F_k : (\bar P_k V', D'_k, \bar V'_k) \cdots \r  (\bar P_k V, D_k, \bar V_k)$$
which commute with the projections, more specifically for all $0 \leq l
\leq k-1$ one has
$$ \pi_{l,k} \circ F_k = F_l \circ \pi_{l,k}'.$$
These maps in turn induce meromorphic maps
$$(F_{k})_* : {\cal O}_{ \bar P_{k+1} V'}(-1) \cdots \r
{\cal O}_{ \bar P_{k+1} V}(-1)  $$
 (holomorphic where $F_k$ is) which also commute with the projections.
\item[b)] If the differential map $F_*: \bar V' \r F^{-1}(\bar V)$ is
injective over a point
$x_0 \in \bar X'$, then there exists a neighborhood $U$ of $x_0$ in $\bar X'$
over which the maps $F_k$ are log-directed morphisms and the induced maps
$$F_k: \bar P_k V' \r F^{-1}(\bar P_k V)$$
are holomorphic embeddings and the induced maps between line bundles
$$(F_{k})_* : {\cal O}_{ \bar P_{k+1} V'}(-1) \r
F^{-1}({\cal O}_{ \bar P_{k+1} V}(-1))  $$
over these embeddings are injective.
\item[c)] If the differential map $F_*: \bar V' \r F^{-1}(\bar V)$ is
bijective in a point
$x_0 \in \bar X'$, then there exists a neighborhood $U$ of $x_0$ in $\bar X'$
over which the maps $F_k$ are log-directed morphisms and the induced maps
$$F_k: \bar P_k V' \r F^{-1}(\bar P_k V)$$
$$(F_{k})_* : {\cal O}_{ \bar P_{k+1} V'}(-1) \r
F^{-1}({\cal O}_{ \bar P_{k+1} V}(-1))  $$
are all bundle isomorphisms over $U$.
\end{itemize}
\end{prop}
Combining  Proposition \ref{functoriality}  with Proposition \ref{pro}
yields the following, which we will use in the
next section
to study $\bar P_kV$  by studying $K^{-1}(\bar P(\cz^r \se E))$:
\begin{prop} \label{localtriv}
Let $(\bar X, D, \bar V)$ be a log-directed manifold.
There exists a neighborhood $U$ of $x_0$ in $\bar X$ and a log-directed
projection
$$K: (\bar X,D,\bar V)|_U \r  (\cz^r, E, \bar T \cz^r) \; ; \; \;
(z_1,..., z_n) \r (z_1,...z_a,z_{l+1},...,z_{l+b}),$$
with $E= \{ z_1...z_a=0 \}$ and $a+b=r= {\rm rank}\, V$, which induces
$$ (K_{k})_*: {\cal O}_{\bar P_{k+1} V}(-1)|_U
 \rightarrow       K^{-1}({\cal O}_{\bar P_{k+1} (\cz^r \se E)}(-1)),$$
$$K_k: \bar P_k V|_U \r K^{-1}(\bar P_k (\cz^r \se E))$$
as bundle isomorphisms.
\hspace*{\fill} $\Box$
\end{prop}
\noindent {\bf Proof for a) } We proceed by induction on $k$. The case
$k=0$ holds by assumption.
Assume the case $k-1$ holds. If
$(F_{k-1})_* : \bar V_{k-1}' \r  \bar V_{k-1}$,  define
$$F_k := \pr ((F_{k-1})_*) : \bar P_kV' \r \bar P_kV.$$
Then by definition of $D_k'= (\pi_k')^{-1} D'_{k-1}$
and $D_k= \pi_k^{-1} D_{k-1}$ the map
$F_k$ is a log-meromorphic morphism which commutes with projections. Let
$(F_k)_* : \bar TP_kV' \r \bar TP_kV$
be the log-differential map defined as in equation (\ref{tt3}).
If $\xi \in ( \bar V_k')_{(w, [v])}$ with $w \in \bar P_{k-1}V'$ and $v \in
(\bar V'_{k-1})_w$, then
$$(\pi_k)_*((F_k)_* \xi ) = (\pi_k \circ F_k)_* \xi =
(F_{k-1} \circ \pi'_{k-1})_* \xi $$
$$= (F_{k-1})_* ( \pi'_{k-1})_* \xi
\in (F_{k-1})_* \cz v = \cz ((F_{k-1})_* v),$$
hence, $(F_k)_* \xi \in ( \bar V_k)_{(F_{k-1}(w), [(F_{k-1})_* v ])}$,
so $F_k$ is a log-directed meromorphic morphism.
The second part of the assertion is clear.\\
{\bf Proof for b and c) }  $(F_*)_x : (\bar V')_x \r
(\bar V)_{F(x)}$
remains injective (respectively bijective)  for all $x$ in a neighborhood
$U(x_0)$.
By the bundle
structures it suffices to prove that for  all $x \in U(x_0)$, the maps
$F_k: (\bar P_kV')_x \r (\bar P_kV)_{F(x)}$ are holomorphic embeddings
(respectively biholomorphic maps) and the maps
$ (F_k)_*: (V_k')_x \r (V_k)_{F(x)}$ are injective (respectively bijective)
bundle maps over them. We prove this by induction on $k$.
  The case $k=0$ holds by assumption.
Assume the case $k-1$ holds. By projectivizing the injective (respectively
bijective) bundle map $(F_{k-1})_*$ and by  a), we get that $F_k : (\bar
P_kV')_x
\r (\bar P_kV)_{F(x)}$ is an injective (respectively bijective) log-directed
morphism. Furthermore, since $F_{k-1} : (\bar P_{k-1}V')_x
\r (\bar P_{k-1}V)_{F(x)}$ is a holomorphic embedding (resp biholomorphic)
and since $(F_{k-1})_* \bar V_k' \s \bar V_k|_{F_{k-1}(\bar P_{k-1}V')}$
is a holomorphic subbundle (respectively the same bundle) the map $F_k$ is an
embedding (respectively biholomorphic). It remains to show that $(F_k)_*$
is again
injective (respectively bijective). Let $\xi \in (\bar V_k')_w$ for $w \in
\bar P_kV'$
 such that $(F_k)_* \xi =0$. Then
$$ 0=(\pi_k)_*  (F_k)_* \xi = (F_{k-1})_*  (\pi'_k)_* \xi.$$
Since the map $(F_{k-1})_*$ is injective  we get $ \xi \in \ker (\pi_k)_*$.
But the subbundle $\ker (\pi_k)_* \s \bar V_k' \s \bar TP_kV'$ can be
canonically identified
with the relative tangent bundle $T_{\bar P_kV'/ \bar P_{k-1}V'}$, which is
a subbundle of $T \bar P_KV'$.
Since we have shown that $F_k$ is a holomorphic embedding, $(F_k)_*$ is
injective
on $(T\bar P_kV')_w$, which contains $\xi$.
As  $(F_k)_* \xi =0$ this forces $\xi =0$. So
 $(F_k)_*$ is injective on $(\bar V_k')_w$. Moreover, if the assumption in
c) holds,
then ${\rm rank}\, V_k' = {\rm rank}\, V'
={\rm rank}\, V ={\rm rank}\, V_k $, and so
$(F_k)_*$ is bijective.
 \qed

\begin{cor} \label{corfunctoriality}
 We have
$F_k((\bar P_k V')^{\rm sing}) \s \bar P_k V^{\rm sing}.$
Moreover, if $F_*: V' \r F^{-1}(V)$ is injective at a
point
$x_0 \in \bar X'$, then there is a neighborhood of $x_0$
over which
$(\bar P_k V')^{\rm sing}$ is isomorphic to $F_k^{-1}(\bar P_k V^{\rm sing})$.
\end{cor}
\noindent {\bf Proof }
By the definition of the singular locus and   Proposition
\ref{functoriality} a)
it suffices to show $F_j (\Gamma_j') \s \Gamma_j$ (respectively
$F_j^{-1} (\Gamma_j) = \Gamma_j'$) for $2 \leq j \leq k$.
Moreover, since the $\Gamma_j'$ and $\Gamma_j$ are divisors without
vertical components,
it suffices to prove the assertions where all maps $F_i$, $i \leq j$ are
holomorphic. The first assertion follows immediately  from the definitions
of $\Gamma_j'$ and $\Gamma_j$ and the equation
\begin{equation} \label{pf3}
(\pi_{l-1})_* \circ (F_{l-1})_* = (F_{l-2})_* \circ (\pi'_{l-1})_*.
\end{equation}
The second follows from this equation and the injectivity
of  $(F_{l-2})_*$.
\qed

\begin{cor} \label{subm}
Let $(\bar X,D, \bar V)$ be a log-directed manifold.
If $\bar V \s \bar W \s \bar TX$ are holomorphic subbundles, then we have
natural
inclusions of submanifolds
$$\bar P_kV \s \bar P_kW \s \bar P_kX$$
and the associated maps  over these inclusions of the line bundles
$${\cal O}_{ \bar P_{k} V}(-1) \s {\cal O}_{ \bar P_{k} W}(-1) \s
{\cal O}_{ \bar P_{k} X}(-1)$$
are line bundle restrictions.
\end{cor}
\noindent {\bf Proof }
We apply Proposition \ref{functoriality} to the log-directed morphism
$i:(\bar X, D, \bar V) \r (\bar X, D, \bar W)$, where $i: \bar X \r \bar X$
is the identity map and induces the bundle inclusion $i_*: \bar V \r \bar W$.
By  Proposition \ref{functoriality} a) and b) we obtain a log-directed morphism
$i_k : \bar P_kV \r \bar P_kW$ which locally over $\bar X$ is, moreover,
a holomorphic embedding $i_k: \bar P_kV \r i^{-1}(\bar P_kW) =
\bar P_kW$. Hence,  $i_k : \bar P_kV \r \bar P_kW$ is a holomorphic embedding.
The other statements follow in a similar way.
 \qed

\section{Log-directed jet differentials}
\setcounter{equation}{0}
\subsection{Demailly-Semple jet bundles and jet differentials}
In this subsection we recall parts of some basic results of
the work \cite{Dem} of Demailly on his construction of the Demailly-Semple
jet bundles, which we
generalize to logarithmic setting in the next subsections.

Let $(X,V)$ be a directed manifold.
Without loss of generality we assume $r = {\rm rank}\, V \geq 2$ in section 3,
for the situation is trivial otherwise.
Let
$$G_k = J_k \cz_0^{\rm reg}= \{t \r \phi (t)=\sum_{i=1}^k a_i t^i\: , \; \;
a_1 \in \cz^*, \; a_i \in \cz , \; i \geq 2 \} $$
be the {\it group of reparametrizations}. Elements $\phi \in G_k$ act on
$J_kV$ as holomorphic automorphisms by
$$\phi: J_kV \r J_kV; \; j_k(f) \r j_k(f \circ \phi).$$
In particular, $\cz^*$ acts on $J_kV$.

Every nonconstant germ
$f:(\cz,0)\r X$ tangent to $V$ lifts to a unique germ
$f_{[k]}:(\cz,0)\r P_kV$ tangent to $V_k$.
$f_{[k]}$ can be defined inductively to be the
projectivization of $f^{\prime}_{[k-1]}:(\cz,0)\r V_{k-1}$.
As such we also have a germ $$f^{\prime}_{[k-1]}:(\cz,0)\r
{\cal O}_{ P_kV}(-1).$$
This construction is actually a special case of our construction in Proposition
 \ref{functoriality}, since $P_k \cz = \cz$.
We get, moreover:
\begin{prop} \label{Fof}
Let $F: (X',V') \r (X,V)$ be a directed morphism. Let $f:(\cz,0)\r X'$
be  a germ tangent to $V'$ such that the germ $F \circ f:(\cz,0)\r X$
is nonconstant. Then there exists a neighborhood $U$ of $0 \in \cz$
such that for all $t \in U$, $t \not= 0$, and for all $k \geq 0$, the map
$F_k : P_kV' \r P_kV$
(see Proposition \ref{functoriality}) is a morphism around $f_{[k]}(t)$,
and we have on $U$:
\begin{equation} \label{pf4}
(F \circ f)_{[k]} = F_k  \circ (f_{[k]}).
\end{equation}
\end{prop}
\noindent {\bf Proof }
Since the germ $F \circ f:(\cz,0)\r X$
is nonconstant, we can find a neighborhood $U$ of $0 \in \cz$
such that  $(F \circ f)'(t) \not= 0$ for all $t \in U$, $t \not= 0$.
{From} the equation $(F \circ f)'(t)= (\pi_{0,k})_*
\circ (F \circ f)_{[k]}'(t)$, we get that $(F \circ f)_{[k]}'(t)\not= 0$ for
all $k \geq 0$.
We now proceed by induction on $k$. The case $k=0$ is trivial. Assume the case
$k-1$. This means that
for all $t \in U$, $t \not= 0$,  the map $F_{k-1} : P_{k-1}V' \r P_{k-1}V$
is a morphism at $f_{[k-1]}(t)$,
and we have on $U$:
$$(F \circ f)_{[k-1]} = F_{k-1} (f_{[k-1]}).$$
Taking the derivative, we obtain
$$(F \circ f)_{[k-1]}'(t) = (F_{k-1})_* (f_{[k-1]})'(t).$$
Now the left hand side is nonzero for $t \not= 0$, so the right hand side
is nonzero, too, and we just can projectivize and obtain the assertion for
$t \not= 0$. Finally equation (\ref{pf4}) still holds for $t=0$ by analytic
continuation.
\qed
\medskip

 The {\it bundle of directed invariant jet differentials of order $k$
and degree $m$}, denoted by $E_{k,m}V^*$, is defined as follows:
Its sheaf of sections ${\cal O}(E_{k,m}V^*)$ over $X$
consists of holomorphic functions on $J_kV|_O$   which satisfy
\begin{equation}\label{invj}
Q(j_k(f\circ\phi))=\phi^\prime(0)^m Q(j_k(f))\ \,\forall\,
j_k(f)\in J_kV^{\rm reg}|_O \mbox{ and }\phi\in G_k
\end{equation}
as $O$ varies over the open subsets of $X$.
We remark that equation (\ref{invj}) implies that the functions $Q$,
restricted to the fibers of $J_kV$, are polynomials of weighted degree
$m$ with respect to the $\cz^*$-action,
so that our definition coincides with the usual one.

\begin{theo}[Demailly (\cite{Dem})] \label{demailly}
Let $(X,V)$ be a directed manifold.
\begin{itemize}
\item[a)] The maps
$$ \tilde{\alpha}_k: J_kV^{\rm reg} \r {\cal O}_{ P_kV}(-1)^{\rm reg}
\: , \; \; j_k(f) \r f^{\prime}_{[k-1]}(0),$$
$$ {\alpha}_k: J_kV^{\rm reg} \r { P_kV}^{\rm reg}
\: , \; \; j_k(f) \r f_{[k]}(0)$$
are well defined, holomorphic and surjective.
\item[b)] If $\phi \in G_k$ is a reparametrization, one has
$$ (f \circ \phi)^{\prime}_{[k-1]}(0)=f^{\prime}_{[k-1]}(0) \cdot
\phi^{\prime}(0), $$
$$ (f \circ \phi)_{[k]}(0)=f_{[k]}(0). $$
\item[c)] The quotient $J_kV^{\rm reg} / G_k$ of $J_kV^{\rm reg}$ by
$G_k$ has the structure of a locally trivial fiber bundle over $X$,
and the map
$$\alpha_k / G_k\: : J_kV^{\rm reg} / G_k \r P_kV$$
is a holomorphic embedding which identifies $J_kV^{\rm reg} / G_k$
with $P_kV^{\rm reg}$.
\item[d)]
The direct image sheaf
$$(\pi_{0,k})_* {\cal O}_{ P_kV}(m) \simeq {\cal O}(E_{k,m}V^*)$$
can be identified with the sheaf ${\cal O}(E_{k,m}V^*)$.
\end{itemize}
\end{theo}

\begin{cor} \label{cordemailly}
$ $\\[-.2in]
\begin{itemize}
\item[1)]
The group $G_k^o =  \{t \r \phi (t)=\sum_{i=1}^k a_i t^i\: , \; \;
a_1=1\} $ acts transitively on the fibers of $\tilde{\alpha}_k$.
\item[2)] The maps $\tilde{\alpha}_k$ and $\alpha_k$ are holomorphic
submersions.
\end{itemize}
\end{cor}
\noindent {\bf Proof for 1) } By  Theorem \ref{demailly} b) and c), the group
$G_k$ acts transitively
on the fibers of $\alpha_k$. So for any two points $p$ and $q$ of a fixed fiber
of $\tilde{\alpha}_k$ we find $\phi \in G_k$ such that  $\phi (p)=q$. Again by
 b) we have
$\phi' (0)=1$, so $\phi \in G_k^o$.\\
{\bf Proof for 2) } It suffices to prove the statement for $\alpha_k$,
since by the action of
$\cz^* = G_k / G_k^o$ on $J_kV^{\rm reg}$ and  Theorem \ref{demailly} b), it
also follows for $\tilde{\alpha}_k$. The assertion is equivalent with the
existence
of local sections for $\alpha_k$ through every point $j_k(f) \in  J_kV^{\rm
reg}$.\\
Let $j_k(f) \in  J_kV^{\rm reg}$ be given, and let $w_0 = \alpha ( j_k(f))$.
Since $f'(0) \not= 0$, we get that $f_{[k-1]}'(0) \not= 0$. Then by Corollary
5.12 in Demailly's  paper \cite{Dem} and its proof we find a neighborhood
$U(w_0) \s P_kV^{\rm reg}$ and a holomorphic family of germs $(f_w)$,
$w \in U(w_0)$,
such that $(f_w)_{[k]}(0)=w$ and $f_{w_0}=f$. After possibly shrinking
$U(w_0)$, we may assume that $f_w'(0) \not= 0$ for all $w \in U(w_0)$.
Thus $w \mapsto j_k(f_w)$ defines a local holomorphic section
$s:U(w_0) \r J_kV^{\rm reg}; \: w \mapsto j_k(f_w)$ with
$s(\alpha_k(j_k(f))) = s(w_0) = j_k(f_{w_0}) = j_k(f)$.
\qed

\subsection{Local trivializations}
\begin{prop} \label{pr32}
Let $z_1,...,z_r$ be the standard coordinates of $\cz^r$, let $a \leq r$,
let $E= \{ z_1...z_a=0 \}$ and
$P=( \overbrace{1,...,1}^a, \overbrace{0,...,0}^{r-a}) \in \cz^r$.
\begin{itemize}
\item[a)]
The trivialization of
$\bar J_k(\cz^r \se E)$ by the forms
$\omega^1= {dz_1\over z_1},...,\omega^a= {dz_a\over z_a},\;
\omega^{a+1}= dz_{a+1},...,$ $\omega^{r}= dz_{r}$,
induce an isomorphism
\begin{equation} \label{pr3211}
\bar J_k (\cz^r \se E) \r J_k(\cz^r)_P \times \cz^r
\end{equation}
which respects regular and singular jets and commutes (outside $E$) with~$G_k$.
\item[b)]
For the log-manifold $(\cz^r,E)$ there exists a line bundle isomorphism
\begin{equation} \label{pr321}
{\cal O}_{\bar P_k(\cz^r \se E)}(-1) \r {\cal O}_{P_k \cz^r}(-1)_P \times \cz^r
\end{equation}
which respects regular and singular jets and such that the diagram
\begin{equation} \label{pr322}
\begin{array}{ccc}

 \bar J_k (\cz^r \se E)^{\rm reg}|_{\cz^r \se E}
& \rightarrow    &   J_k (\cz^r)_P^{\rm reg} \times ( \cz^r \se E)      \\
                    &                &             \\
\downarrow\tilde \alpha_k&  &\ \ \ \downarrow (\tilde \alpha_k)_P \times id \\
                    &                &             \\
 {\cal O}_{\bar P_k (\cz^r \se E)}(-1)^{\rm reg}|_{\cz^r \se E}     &
\rightarrow    &
   {\cal O}_{P_k (\cz^r)}(-1)_P^{\rm reg} \times ( \cz^r \se E) \\

\end{array}
\end{equation}
commutes.
\end{itemize}
\end{prop}
\noindent {\bf Proof for a) }
The composition
$$(J_k\cz^r)_P = J_k(\cz^r \se E)_P \hookrightarrow \bar J_k (\cz^r \se E)
\r (\cz^k)^r$$
is an isomorphism, where the last morphism is given by the first factor of
the trivialization
map in equation (\ref{t6}). We compose the isomorphism
of equation (\ref{t6}) with the inverse of the above to obtain the isomorphism
$$\bar J_k (\cz^r \se E) \r J_k(\cz^r)_P \times \cz^r\, ; $$ $$
( (Z_j^i)_{i=1,...,r;j=1,...,k}; x) \r
(( (Z_j^i)_{i=1,...,r; j=1,...,k};P);x)\, .$$
This isomorphism respects regular and singular jets,
since the subset of the singular jets is given in every fiber  by
$\{ Z^i_1=0\, , \: i=1,...,r \}$ by Proposition~\ref{p}.
\smallskip

Let us understand this isomorphism, restricted to $\cz^r \se E$, in a more
geometric way. As in the proof of Proposition \ref{barJ_kV}, let
$$\Psi : \cz^r \r \cz^r; \; (w_1,...,w_r) \r (e^{w_1},...,e^{w_a},
w_{a+1},...,w_r)
=(z_1,...,z_r),$$
 and let $W_j^i$, $i=1,...,r; \; j=1,...,k$ be the components  of
the first part of
$$  (\tilde{dw_1},..., \tilde{dw_r}) \times \pi: J_k (\cz^r) \r
 (\cz^k)^r \times \cz^r.$$
Then we claim that the above isomorphism, restricted to $\cz^r \se E$, is
given by
$$J_k (\cz^r \se E) \r J_k(\cz^r)_P \times \cz^r \se E\, ; \;
j_k(f) \r (j_k(\Psi (\Psi^{-1} \circ f  - \Psi^{-1} \circ f(0))), f(0)),$$
where subtraction means subtraction in $\cz^r$. Note that $\Psi^{-1}$ is only
defined up to addition of summands $2 \pi i m, \: m \in \gz$ for the first
$a$ components, but the germ
$\Psi^{-1} \circ f  - \Psi^{-1} \circ f(0)$ is well defined. In fact, by
Lemma \ref{lemW} the above isomorphism is given by trivial shift with respect
to the coordinates $W^i_j$, but this is, by the definition of these
coordinates, just subtraction of the value of the germ for $t=0$.
It follows that the above isomorphism
commutes with the action of $G_k$. In fact, reparametrization does not
depend on
the coordinates and so it commutes with $\Psi$ and $\Psi^{-1}$. Furthermore,
it commutes with subtraction of constants in $\cz^r$. This proves
a).\\
{\bf Proof for b) } We use the following strategy: Using some results
of Demailly's paper \cite{Dem} we first define the isomorphism of equation
(\ref{pr321})
on $\cz^r \se E$ similarly to our geometric way above. It is then easy
to verify the diagram of equation (\ref{pr322}). We then extend this
isomorphism over the complement of a divisor in the bundle which does not
contain
any entire fiber over  $\cz^r$. For this,
we introduce an explicit local coordinate chart in
$\bar P_k (\cz^r \setminus E)$ the complement of which is a divisor which
does not contain any fiber over
$\cz^r$. In order to extend
over the remaining codimension two locus we use the fact that our objects
are defined inductively by projectivizing
vector bundles, and that for vector bundle maps, the Riemann
Extension Theorem holds.
This way is not the shortest possible
(see Lemma \ref{sat}, which gives
a much shorter and intrinsic proof of this extension over the
divisor $E$), but  it
explains well the geometric contents of the isomorphism in equation
(\ref{pr321}) via explicit local coordinates
(see also Corollary \ref{coor}), which is
useful for applications.

By Corollary 5.12 of \cite{Dem}, for all points $w \in
P_k( \cz^r \se E)$,
there exists a germ $f:(\cz,0) \r \cz^r \se E$ such that $f_{[k]}(0)=w$ and
$f_{[k-1]}'(0)\not=0$. We claim that by composing this germ with the map
$t\r at + t^2$,
$a \in \cz$, the vector $f_{[k-1]}'(0)$ can be made equal to an  arbitrary
vector in the complex line
${\cal O}_{ P_k (\cz^r \se E)}(-1)$ over $w$. This follows from
Theorem \ref{demailly}  b) for $a \not= 0$. Since  the image of the germ
$f_{[k]}$ does not change, we get, after an easy computation,
 that  $f_{[k-1]}'(0)=0$ for $a=0$. So every vector of
$ {\cal O}_{ P_k (\cz^r \se E)}(-1)$ is obtained this way.

As above, the map
$$ {\cal O}_{ P_k (\cz^r \se E)}(-1)\r ({\cal O}_{ P_k (\cz^r \se E)}(-1))_P
 \times (\cz^r \se E)\, ; \;$$ $$
f_{[k-1]}'(0) \r (\Psi (\Psi^{-1} \circ f  - \Psi^{-1} \circ f(0))_{[k-1]}'(0);
 f(0))$$
is a well defined isomorphism, its inverse being
given by the well defined map
$$ ({\cal O}_{ P_k (\cz^r \se E)}(-1))_P
 \times \cz^r \se E \r {\cal O}_{ P_k (\cz^r \se E)}(-1)
\, ; \;$$ $$
(f_{[k-1]}'(0),p) \r \Psi (\Psi^{-1} \circ f  + \Psi^{-1}(p))_{[k-1]}'(0)\; .$$
Now, by Proposition 5.11 of  \cite{Dem},
the singular locus of
$ {\cal O}_{ P_k (\cz^r \se E)}(-1)$ can be characterized by
$f_{[k-1]}'=0$ along with the multiplicities
 of $f_{[j]}$, $j=0,...,k-1$, which remain invariant under changes of
coordinates
or additions by constants. So this isomorphism respects the regular and
singular jet loci. So we can restrict it to the regular loci. If now
$j_k(f) \in J_k (\cz^r \se E)^{\rm reg}$, then this is mapped to
$(\Psi (\Psi^{-1} \circ f  - \Psi^{-1} \circ f(0))_{[k-1]}'(0),f(0))$
by both compositions of the maps in the diagram in equation (\ref{pr322}),
so this diagram commutes. \\

We now carry out the above strategy via the following lemma.
Let $\bar V_k^{\rm reg} = \bar V_k \se  \bar P_k V|_{\bar P_kV^{\rm reg}}$,
where the $\bar P_k V$ denotes the zero section of $\bar V_k$.
Then we have a canonical identification
$${\cal O}_{ \bar P_k V}(-1)^{\rm reg}
\stackrel{\sim}{\longrightarrow} \bar V_{k-1}^{\rm
reg}\, .$$
We now introduce $r$ coordinate
charts on $\bar V_{k-1}^{\rm reg}$
in the same way as Demailly did for
$P_kV^{\rm reg}$ in equations (4.9) and (5.7) and Theorem 6.8 of \cite{Dem}.

\begin{lem} \label{chart}
Let
$( (Z_j^i)_{i=1,...,r;j=1,...,k}; (z_1,...,z_r))$ be the coordinates of
$\bar J_k(\cz \se E)$, and let
$$\bar A_{k,r} =
(\bigcap_{j=2}^k \{ Z^r_j =0 \})
\cap ( \bar J_k(\cz \se E) \se \{Z_1^r=0
\})\, .$$
Then the map
\begin{equation} \label{ch1}
\tilde \alpha_k : \bar A_{k,r}|_{\cz^r \se E} \r \bar V_{k-1}|_{\cz^r \se E}
\end{equation}
extends over $E$ to a map which is biholomorphic onto its image
$\tilde \alpha_k (A_{k,r})$, such that
this image contains the complement of a
divisor in $\bar V_{k-1}$ which is nowhere dense in  $(\bar
V_{k-1})_x$ for all $x \in \cz^r$.
More precisely:\\

\noindent {\bf Claim S(k):}
$\bar V_{k-1} \r \bar P_{k-1}(\cz^r \se E)$
is a vector bundle of rank $r$ over a
$(k-1)$-stage tower of $\pr^{r-1}$-bundles, and we can introduce
inhomogenous coordinates  on these bundles corresponding to the coordinates
$(z_1,...,z_r)$ of $\cz^r$,
in which the map $\tilde \alpha_k$ of equation (\ref{ch1})  is given by
$$( (Z_j^i)_{i=1,...,r-1;j=1,...,k},Z_1^r; (z_1,...,z_r)) \r
((\frac{Z^1_1}{Z_1^r},...,\frac{Z^{r-1}_1}{Z_1^r}),
(\frac{Z^1_2}{(Z_1^r)^2},...,\frac{Z^{r-1}_2}{(Z_1^r)^2}),...,$$ $$
(\frac{Z^1_{k-1}}{(Z_1^r)^{k-1}},...,\frac{Z^{r-1}_{k-1}}{(Z_1^r)^{k-1}}),
(\frac{Z^1_{k}}{(Z_1^r)^{k-1}},...,\frac{Z^{r-1}_{k}}{(Z_1^r)^{k-1}},
Z_1^r);(z_1,...,z_r))\, .$$
\end{lem}
{\bf Proof of Lemma \ref{chart} }
It suffices to prove {\bf S(k)}.
We prove this by induction. The statement
 {\bf S(1)} is trivial. Assume by
induction that {\bf S(k-1)} holds.
Then the corresponding inhomogenous
coordinates of $\bar P_{k-1}(\cz^r \se E)$ are given by
\begin{equation} \label{ch2}
((\frac{Z^1_1}{Z_1^r},...,\frac{Z^{r-1}_1}{Z_1^r}),
(\frac{Z^1_2}{(Z_1^r)^2},...,\frac{Z^{r-1}_2}{(Z_1^r)^2}),...,
(\frac{Z^1_{k-1}}{(Z_1^r)^{k-1}},...,\frac{Z^{r-1}_{k-1}}{(Z_1^r)^{k-1}});(z_1,.
..,z_r))\, .
\end{equation}
In order to find coordinates over this
affine chart, we proceed in two steps:

The first step is to find the coordinates of
the logarithmic tangent bundle
$$\bar T P_{k-1}(\cz^r \se E) \r
\bar  P_{k-1}(\cz^r \se E)$$ over our
affine chart. Note  that, in the coordinates of equation (\ref{ch2}), the
divisor $E_{k-1}= \pi_{0,k-1}(E)$ is
given by $\{z_1  ...  z_a =0 \}$ and, hence, is independent of
any of the fiber coordinates $Z_j^i$ or of their quotients
${Z^i_j}/{(Z^r_1)^j}$. So the
 coordinates of of $\bar T P_{k-1}(\cz^r \se E)$ are given by those
of equation (\ref{ch2}) and their
 differentials, except for
$z_1,...,z_a$, where
the log-differentials are needed.

The second step is to restrict this coordinate system to the
subbundle $\bar V_{k-1} \s \bar T P_{k-1}(\cz^r \se E)$. By the definition
of this subbundle in equation (\ref{O1})
(see also Demailly's equation (5.7) in
\cite{Dem}) we choose the differentials of the $r-1$ coordinate functions
of equation (\ref{ch2}) which
describe the fibers of the map
$\pi_{k-1}: \bar P_{k-1}(\cz^r \se E)
\r \bar P_{k-2}(\cz^r \se E)$ and of
an extra component which corresponds to how we have chosen the inhomogenous
coordinates: These are the
(nonlog-) differentials
$d(\frac{Z_{k-1}^i}{(Z_1^r)^{k-1}})$,
$i=1,...,r$ plus the extra component, corresponding to the log- (in case
$a=r$) or nonlog- (in case $a < r$) differential of $z_r$, which is
$Z_1^r$ (see equations (\ref{t1}) and
(\ref{t6})).

It remains to express these coordinates without the differential
as in claim {\bf S(k)}.
We have for $1 \leq i \leq r-1$:
$$ d(\frac{Z_{k-1}^i}{(Z_1^r)^{k-1}})=
\frac{dZ_{k-1}^i}{(Z_1^r)^{k-1}}-
\frac{Z_{k-1}^i}{(Z_1^r)^{k-1}} \cdot
(k-1) \frac{dZ_{1}^r}{Z_1^r}\, .$$
By equations (\ref{t1}) and (\ref{t6}) we get
$$dZ_{k-1}^i(j_k(f)) = (\frac{d}{dt}
\frac{d^{k-2}}{dt^{k-2}}
 \frac{f^* \omega^i}{dt})|_{t=0}
= (\frac{d^{k-1}}{dt^{k-1}}
 \frac{f^* \omega^i}{dt})|_{t=0}=
Z_k^i(j_k(f))
\, .$$
As we only work on the
submanifold $\bar A_{k,r}$,
we have $j_k(z_r \circ f)= a_0(x)+a_1(x)t$.
We now again use equations (\ref{t1}) and (\ref{t6}),
and distinguish two cases:\\
If $a <r$, then
$$dZ_{r}^1(j_k(f)) = (\frac{d}{dt}
 \frac{f^* dz_r}{dt})|_{t=0}
= (\frac{d}{dt}a_1(x))|_{t=0}= 0\, .$$
If $a=r$, then
$$dZ_{r}^1(j_k(f)) = (\frac{d}{dt}
 \frac{f^* dz_r}{(z_r \circ f)dt})|_{t=0}
= (\frac{d}{dt}\frac{a_1(x)}{a_0(x)+a_1(x)t})|_{t=0}=$$
$$=-(\frac{a_1(x)}{a_0(x)})^2 =
-( \frac{f^* dz_r}{(z_r \circ f)dt}|_{t=0})^2
= -(Z_1^r)^2(j_k(f))\, .$$
So we have
$$d(\frac{Z_{k-1}^i}{(Z_1^r)^{k-1}})=
\frac{Z_{k}^i}{(Z_1^r)^{k-1}}+ (k-1) \cdot
\left\{
\begin{array}{ccc}
0 & : & a<r\\
Z_1^r & : & a=r
\end{array} \right. \: .$$
These coordinates can be expressed in those of claim {\bf S(k)} and vice
versa. \qed
\medskip
Now we use Lemma \ref{chart} to extend the diagram in equation
(\ref{pr321}), and,
thus, the isomorphism of equation (\ref{pr321}): The diagram becomes
$$
\begin{array}{ccc}

 \bar A_{k,r}|_{\cz^r \se E}
& \rightarrow    &   (A_{k,r})_P \times ( \cz^r \se E)      \\
                    &                &             \\
\downarrow\tilde \alpha_k&  &\ \ \ \downarrow (\tilde \alpha_k)_P \times id \\
                    &                &             \\
 \bar V_{k-1}|_{\cz^r \se E}     & \rightarrow    &
   (V_{k-1})_P \times ( \cz^r \se E) \\

\end{array}
$$
where the two vertical arrows are now biholomorphic onto their images.
By Lemma \ref{chart} these isomorphisms
extend to isomorphisms over $\cz^r$.
So the isomorphism of equation
(\ref{pr321}) extends over $E$ outside a horizontal
divisor which is nowhere dense in all fibers, giving an isomorphism outside
an analytic set of codimension at least two.

We finally prove, by induction over $k$,
that these isomorphisms extend to
$$
\begin{array}{ccc}

 \bar V_{k-1}
& \rightarrow    &   (V_{k-1})_P \times  \cz^r       \\
                    &                &             \\
\downarrow&  &\ \ \ \downarrow  \\
                    &                &             \\
 \bar P_{k-1}V    & \rightarrow    &
   (P_{k-1}V)_P \times  \cz^r  \\

\end{array}
$$
which induce the desired isomorphisms
of equation (\ref{pr321}). The case
{\bf S(1)} is trivial (there is nothing to extend any more in this case).
Assume by induction that {\bf S(k-1)} is
true. Then by projectivizing we have an
isomorphism
$\bar P_{k-1}V     \rightarrow
   (P_{k-1}V)_P \times  \cz^r$,
and over this we have an isomorphism
$ \bar V_{k-1}     \rightarrow
   (V_{k-1})_P \times  \cz^r$
up to a subvariety of codimension two.
Now this isomorphism extends, since
for vector bundle maps, the Riemann
Extension Theorem holds. (For any point $w \in \bar P_{k-1}V$, take
a dual basis of $\bar V_{k-1}$ around $w$. Then the extension of the maps, in
both directions, is reduced to extension of holomorphic functions once we
compose
these maps with the dual vectors.)
The fact that the extended maps are still inverses to each other follows
from the Identity Theorem. This ends the proof of Proposition~
\ref{pr32}.\qed

\noindent {\bf Important Remark 1 } The local isomorphisms
of equations (\ref{pr3211}) and (\ref{pr321}) are fiber bundle
isomorphisms. But they are
{\it not} induced by (directed) morphisms.
As a result, these local
isomorphisms have a priori no  functoriality, and every compatibility
which one needs has to be proved explicitly, which we proceed to do.

\begin{prop} \label{pr33}
Let $(\bar X, D, \bar V)$ be a log-directed manifold.
Let $x_0 \in \bar X$ and let $U$ and the log-directed projection
$$K: (\bar X,D,\bar V)|_U \r  (\cz^r, E, \bar T \cz^r) \; ; \; \;
(z_1,..., z_n) \r (z_1,...z_a,z_{l+1},...,z_{l+b}),$$
with $E= \{ z_1...z_a=0 \}$ and $a+b=r= {\rm rank}\, V$,
be as in Proposition \ref{localtriv}.
Let, without loss of generality,
 $P=( \overbrace{1,...,1}^l, \overbrace{0,...,0}^{n-l}) \in U$.
\begin{itemize}
\item[a)]
The isomorphisms of Proposition \ref{pr32},
Proposition \ref{barJ_kV} and Proposition \ref{localtriv}
induce isomorphisms
$$\bar J_kV|_U     \rightarrow  J_kV_P \times U $$
and
$${\cal O}_{\bar P_kV}(-1)|_U     \rightarrow
   {\cal O}_{P_kV}(-1)_P \times U$$
respecting regular and singular jets
in such a way that the first isomorphism commutes (outside $D$) with the
action of $G_k$
and that the diagrams
 {\small\begin{equation} \label{pr331}
\begin{array}{rccr}

 K^{-1}(\bar J_k (\cz^r \se E))|_U
& \rightarrow    &   (K^{-1}(J_k \cz^r))_P \times U &     \\
                    &        &        &            \\
   \nwarrow & &   \hfill \nwarrow & \\
                    &            &    &           \\
  & \bar J_kV|_U    & \rightarrow  &
   J_kV_P \times U \\

\end{array}
\end{equation}}

and
{\small \begin{equation} \label{pr332}
 \begin{array}{rccr}

 K^{-1}({\cal O}_{\bar P_k (\cz^r \se E)}(-1))|_U &
 \rightarrow    &   (K^{-1}({\cal O}_{P_k \cz^r}(-1)))_P \times U &      \\
                    &                &       &      \\
   \nwarrow & & \hfill  \nwarrow & \\
                    &                &       &      \\
 &  {\cal O}_{\bar P_kV}(-1)|_U    & \rightarrow  &
   {\cal O}_{P_kV}(-1)_P \times U \\

\end{array}
\end{equation}}
commute.

\item[b)]
Moreover, outside the divisor $D$, they induce the following cubic
diagram
 {\footnotesize
\begin{equation} \label{pr333}
\hspace{-0.5 in}
\begin{array}{cccc}

 K^{-1}(\bar J_k (\cz^r \se E)^{\rm reg})|_{U \se D} &
 \rightarrow    &   (K^{-1}(J_k \cz^r)^{\rm reg})_P \times (U \se D) &      \\
                    &                &       &      \\
   \hfill \nwarrow & & \hfill  \nwarrow & \\
                    &                &       &      \\
  & \bar J_kV^{\rm reg}|_{U \se D}    & \rightarrow  &
   J_kV_P^{\rm reg} \times (U \se D) \\
&&&\\
\downarrow \tilde{\alpha}_k&&\downarrow (\tilde{\alpha}_k)_P \times id&\\
&&&\\
&&&\\
&\downarrow \tilde{\alpha}_k&&\downarrow (\tilde{\alpha}_k)_P \times id\\
&&&\\
&&&\\
 K^{-1}({\cal O}_{\bar P_k (\cz^r \se E)}(-1)^{\rm reg})|_{U \se D} &
 \rightarrow    &   (K^{-1}({\cal O}_{P_k \cz^r}(-1))^{\rm reg})_P \times
(U \se D) &      \\
                    &                &       &      \\
  \hfill \nwarrow & & \hfill  \nwarrow & \\
                    &                &       &      \\
 &  {\cal O}_{\bar P_kV}(-1)^{\rm reg}|_{U \se D}    & \rightarrow  &
   {\cal O}_{P_kV}(-1)_P^{\rm reg} \times (U \se D) \\

\end{array}
\end{equation}}

\item[c)] By combining with the canonical line bundle projections
we get the same isomorphisms and diagrams with
$\bar P_k (\cz^r \se E)$,
$\bar P_kV$ and
$\alpha_k$ instead of
${\cal O}_{\bar P_k (\cz^r \se E)}(-1)$,
${\cal O}_{\bar P_kV}(-1)$ and
$\tilde{\alpha}_k$.
\end{itemize}
\end{prop}
\noindent {\bf Proof for a) } We {\it define} the isomorphisms
$\bar J_kV|_U \r J_kV_P \times U$
respectively
${\cal O}_{\bar P_kV}(-1)|_U     \rightarrow
   {\cal O}_{P_kV}(-1)_P \times U $
by the other three arrows of the respective diagrams.
 In this way we obtain trivializations which,
by definition, make the diagrams commutative. By Proposition \ref{pr32},
Proposition \ref{p}, Proposition \ref{barJ_kV}
and Corollary \ref{corfunctoriality}, the regular
and singular loci are preserved.
The first isomorphism commutes with the action of $G_k$. In fact,
by Proposition \ref{pr32} a) this is true for the upper line
of the diagram in equation (\ref{pr331}). Furthermore, the isomorphism
$K_k$ in the vertical arrows is, outside $D$, just given by
$K_k(j_k(f))=j_k(K \circ f)$, and this trivially commutes with the action
of $G_k$.\\
{\bf Proof for b) } The back side of this cubic diagram (the side with the
$K^{-1}$)
is just the pull back the diagram in equation (\ref{pr322}).
The upper and the lower sides of the cubic diagram are the restrictions
of the diagrams in equations (\ref{pr331}) respectively (\ref{pr332}) to the
regular locus over $U \se D$.  The two vertical arrows
on the front side are  defined by the left hand side respectively
the right hand side of the cubic diagram, so these two sides commute
by definition. It is an easy exercise to see that then the the front side
of the cubic diagram commutes also and, furthermore, that the whole diagram
commutes.\\
{\bf Proof for c) } This is clear from the diagrams.
\qed

\noindent {\bf Important Remark 2 } For all local isomorphisms given
by the horizontal left-to-right arrows in the above diagrams, our Important
Remark 1
 also applies.
However, the local isomorphisms induced by
$K$ are functorial.
\newpage

\begin{cor} \label{coor} $ $\\[-.3in]
\begin{itemize}
\item[a)] The fiber bundles $\bar P_kV$,
${\cal O}_{\bar P_kV}(-1)$ and $\bar V_k$ and their regular and singular
jet loci are all locally trivialized over $\bar X$ in a way
which is compatible, through the maps
$\alpha_k$ respectively $\tilde{\alpha_k}$, with the trivialization of
$\bar J_kV$ by using local
logarithmic coordinates .
\item[b)] Let $U \s \bar X$ and $K$ be like in
Proposition \ref{localtriv}.
Let $\bar A_{k,i} \s \bar J_k(\cz^r \se E)$, $i=1,...,r$, be like in Lemma
\ref{chart}, and let $\bar B_{k,i}
= \bar A_{k,i} \cap \{Z_1^i =1 \}$.
Then there exist $r$ coordinate charts
$$K^{-1}(\bar A_{k,i}) \r \bar V_{k-1}\:
({\rm respectively}\:
K^{-1}(\bar A_{k,i}) \r
{\cal O}_{\bar P_kV}(-1));$$
$$( (Z_j^i)_{i=1,...,r-1;j=1,...,k},Z_1^r; (z_1,...,z_n)) \r
((\frac{Z^1_1}{Z_1^r},...,\frac{Z^{r-1}_1}{Z_1^r}),
(\frac{Z^1_2}{(Z_1^r)^2},...,\frac{Z^{r-1}_2}{(Z_1^r)^2}),...,$$ $$
(\frac{Z^1_{k-1}}{(Z_1^r)^{k-1}},...,\frac{Z^{r-1}_{k-1}}{(Z_1^r)^{k-1}}),
(\frac{Z^1_{k}}{(Z_1^r)^{k-1}},...,\frac{Z^{r-1}_{k}}{(Z_1^r)^{k-1}},
Z_1^r);(z_1,...,z_n))$$
which cover $\bar V_{k-1}^{\rm reg}$
(respectively ${\cal O}_{\bar P_kV}(-1)^{\rm reg}$), and $r$ coordinate charts
$$K^{-1}(\bar B_{k,i}) \r \bar P_{k-1}V;$$
$$( (Z_j^i)_{i=1,...,r-1;j=1,...,k}; (z_1,...,z_n)) \r $$ $$
((Z^1_1,...,Z^{r-1}_1),
(Z^1_2,...,Z^{r-1}_2),...,
(Z^1_{k},...,Z^{r-1}_{k});(z_1,...,z_n))$$
which cover $\bar P_{k-1}V^{\rm reg}$.
\end{itemize}
\end{cor}
{\bf Proof }  a) is contained in Proposition \ref{pr33}. The existence of
the coordinate charts in  b) follows
from Proposition \ref{pr33} and Lemma
\ref{chart}. These charts cover the locus of regular jets by Theorem 3.2, a)
outside $D$. Thus by our local trivializations which are compatible with
the charts
and with the locus of regular jets, these charts cover
the locus of regular jets everywhere. \qed

\noindent {\bf Remark }
The coordinates can also be obtained directly
without Lemma \ref{chart}.

\subsection{Log-directed jets and log-Demailly-Semple jets}
This subsection extends Theorem \ref{demailly}
a), b) and c) to the log-directed case.
\begin{prop} \label{demaillytheo}
Let $(\bar X,D, \bar V)$ be a log-directed manifold.
\begin{itemize}
\item[a)] The maps $\tilde{\alpha}_k$ and $\alpha_k$ of Theorem
\ref{demailly} a)
extend to holomorphic and surjective maps
$$ \tilde{\alpha}_k: \bar J_kV^{\rm reg} \r {\cal O}_{ \bar P_kV}(-1)^{\rm reg}
\: , $$
$$ {\alpha}_k: \bar J_kV^{\rm reg} \r { \bar P_kV}^{\rm reg}\:. $$
\item[b)] The action of $\phi \in G_k$ extends to an automorphism
of $\bar J_kV$ leaving $\bar J_kV^{\rm reg}$ and $\bar J_kV^{\rm sing}$
invariant and  satisfying
$$\tilde{\alpha}_k \circ \phi = \tilde{\alpha} \cdot \phi'(0)\: ,
 \; \alpha_k \circ \phi = \alpha_k\: .$$
\item[c)] The quotient $\bar J_kV^{\rm reg} / G_k$ has the structure of a
locally trivial fiber bundle over $\bar X$,
and the map
$$\alpha_k / G_k\: : \bar J_kV^{\rm reg} / G_k \r  \bar P_kV$$
is a holomorphic embedding which identifies $\bar J_kV^{\rm reg} / G_k$
with $\bar P_kV^{\rm reg}$.
\end{itemize}
\end{prop}
{\bf Proof for a) }
By Theorem \ref{demailly}, the map $\tilde{\alpha}_k$
is defined outside $D$:
\begin{equation} \label{334}
 \tilde{\alpha}_k: \bar J_kV^{\rm reg}|_{\bar X \se D} \r
{\cal O}_{ \bar P_kV}(-1)^{\rm reg}|_{\bar X \se D}.
\end{equation}
Let $x \in D$. By Proposition \ref{pr33} b),
 there exists a neighborhood $U$
of $x$ with
\begin{equation}
\begin{array}{ccc}

  \bar J_kV^{\rm reg}|_{U \se D}
& \rightarrow    &   J_kV_P^{\rm reg} \times (U \se D)      \\
                    &                &             \\
\downarrow\tilde{\alpha}_k&  &\ \ \ \downarrow (\tilde{\alpha}_k)_P \times
id \\
                    &                &             \\
   {\cal O}_{\bar P_kV}(-1)^{\rm reg}|_{U \se D}   & \rightarrow    &
   {\cal O}_{P_kV}(-1)_P^{\rm reg} \times (U \se D)\\
\end{array}
\end{equation}
Here the horizontal arrows are isomorphisms, which, by Proposition
\ref{pr33}a),
extend as isomorphisms over $U$, and $(\tilde{\alpha}_k)_P \times id$ is
clearly extendable to $U$ to a surjective holomorphic map
 on the right hand side. So $\tilde{\alpha_k}$ is also
extendable to a surjective holomorphic map over
$U$ on the left hand side. Since $x \in D$ is arbitrary,
and since by equation (\ref{334}) the extension of $\tilde{\alpha}_k$
is unique if it exists, we obtain a well defined surjective holomorphic
map
\begin{equation} \label{335}
 \tilde{\alpha}_k: \bar J_kV^{\rm reg} \r
{\cal O}_{ \bar P_kV}(-1)^{\rm reg}.
\end{equation}
By combining with the canonical line bundle projections
we get in the same way a surjective holomorphic map
\begin{equation} \label{336}
 \alpha_k: \bar J_kV^{\rm reg} \r
\bar P_kV^{\rm reg}
\end{equation}
which extends the corresponding map $\alpha_k$ of
Theorem \ref{demailly} from $\bar X \se D$ to $\bar X$.\\
{\bf Proof for b) }
If $\phi \in G_k$ is a reparametrization, one has on
$\bar J_kV^{\rm reg}|_{\bar X \se D}$ by Theorem \ref{demailly}:
 \begin{equation} \label{336a}
\tilde{\alpha}_k \circ \phi=\tilde{\alpha}_k \cdot
\phi^{\prime}(0) \; ,\; \;
 \alpha_k \circ \phi = \alpha_k,
\end{equation}
where in the first equation the multiplication
$\tilde{\alpha_k} \cdot
\phi^{\prime}(0)$ denotes the multiplication with scalars in the line
bundle ${\cal O}_{ \bar P_kV}(-1)^{\rm reg}|_{\bar X \se D}$.
By Proposition \ref{pr33} a), we have the  diagram
\begin{equation} \label{337}
\begin{array}{ccc}

  \bar J_kV|_{U \se D}
& \rightarrow    &   J_kV_P \times (U \se D)      \\
                    &                &             \\
\downarrow  \phi&  &\ \ \ \downarrow  \phi_P \times id \\
                    &                &             \\
\bar J_kV|_{U \se D}
& \rightarrow    &   J_kV_P \times (U \se D)      \\
\end{array}.
\end{equation}
By a similar argument as in  a), the map $ \phi$ extends
to a holomorphic automorphism on $\bar J_kV$.
{From} this diagram, it also follows
that $\phi$ maps $\bar J_kV^{\rm reg}$ onto itself, since all
arrows of this diagram
preserve regular and singular jets, and by Proposition \ref{p} this
remains true over $D$.
Finally, equation (\ref{336a}) extends from $\bar J_kV^{\rm reg}|_{\bar X
\se D}$ to
$\bar J_kV^{\rm reg}$ by the Identity Theorem.\\
{\bf Proof for c) }
  By  b), the quotient
$\bar J_kV^{\rm reg}/G_k$ is well defined (as set).
By the diagrams of equations (\ref{336a}) and
(\ref{337}), we obtain from  Proposition \ref{pr33} c):
\begin{equation} \label{338}
\begin{array}{ccc}

  \bar J_kV^{\rm reg}/G_k|_{U \se D}
& \rightarrow    &   (J_kV^{\rm reg}/G_k)_P \times (U \se D)      \\
                    &                &             \\
\downarrow\alpha_k/G_k&  &\ \ \ \downarrow (\alpha_k/G_k)_P \times id \\
                    &                &             \\
   \bar P_kV^{\rm reg}|_{U \se D}   & \rightarrow    &
   P_kV_P^{\rm reg} \times (U \se D)\\
\end{array}
\end{equation}
 By Demailly (\cite{Dem}), the vertical arrows in this diagram are
isomorphisms. By a similar argument as in  a), one obtains
a holomorphic isomorphism
\begin{equation} \label{339}
\alpha_k/G_k : \bar J_kV^{\rm reg}/G_k \r \bar P_kV^{\rm reg}
\end{equation}
over $\bar X$.  Equation (\ref{338}) shows that this isomorphism
makes $\bar J_kV^{\rm reg}/G_k$ into a holomorphic fiber bundle over $\bar X$.
\qed

\subsection{Characterization of log-directed jet differentials}
In this subsection we generalize  Theorem \ref{demailly} d).
More precisely we prove:
\begin{prop} \label{dlog}  A holomorphic (respectively meromorphic)
function $Q$ on $\bar J_k V|_{ O}$ for some connected open subset
 $ O \s \bar X$  which satisfies
\begin{equation} \label{dlog1}
Q(j_k(f\circ\phi))=\phi^\prime(0)^m Q(j_k(f))\ \,\forall\,
j_k(f)\in J_kV^{\rm reg}\mbox{ and } \forall \phi\in G_k
\end{equation}
 over
some open subset of $ O'$ of $ O \se D$  defines a holomorphic
(respectively meromorphic)
section of ${\cal O}_{\bar P_kV}(m)$ over $ O$,
and vice versa.
\end{prop}
\noindent {\bf Proof }
Let $Q : \bar J_kV|_O \r \cz$ be a meromorphic function
which satisfies
\begin{equation} \label{pf5}
Q \circ \phi = \phi^\prime(0)^m Q \;\; \forall \phi\in G_k
\end{equation} over $O'$.
Since $\bar J_kV^{\rm reg}|_O$ is connected, equation (\ref{pf5})
holds over $O$ by the Identity Theorem.
Since $\tilde \alpha_k$ and $\alpha_k$ were obtained over $D$
by trivial extensions in the diagrams of Proposition \ref{pr33}, the
results of Corollary \ref{cordemailly} extend also over $D$.
In particular, $G_k^o$ of Corollary \ref{cordemailly}
acts transitively on the fibers of $\tilde
\alpha_k$.
Since the function $Q$ is invariant under the
action of $G_k^o$  by equation (\ref{pf5}),
there exists a Zariski-densely defined function
$\tilde Q :{\cal O}_{\bar P_kV}(-1)^{\rm reg}|_{O} \r \cz$ such that
$Q = \tilde Q \circ \tilde \alpha_k$. Again by Corollary \ref{cordemailly},
$\tilde \alpha_k$ has local holomorphic sections everywhere. So
$\tilde Q$ is a meromorphic function on ${\cal O}_{\bar P_kV}(-1)^{\rm
reg}|_{O}$.
By equation (\ref{pf5}), this function is $m$-linear with respect to
the $\cz^*$-action and so corresponds to a meromorphic section $s$ of
${\cal O}_{\bar P_kV}(m)^{\rm reg}|_{U}$.
In order to extend this section to the singular locus, we have to redo
an argument of Demailly (\cite{Dem}):
In a neighborhood $W$ of any point $w_0 \in \bar P_kV|_{O \se D}$ we can
find a holomorphic
family of germs $f_w$ such that $(f_w)_{[k]}(0)=w$ and $(f_w)'_{[k-1]}(0)
\not= 0$.
Then we get
$s(w)=Q(f_w',...,f_w^{(k)})(0)((f_w)'_{[k-1]})^m$
on $W \cap \bar P_kV|_{O \se D}$.
Now the right hand side extends to a section of
${\cal O}_{\bar P_kV}(m)|_{U \se D}$, so the left hand side does, too.
So $s$ is a meromorphic section of
${\cal O}_{\bar P_kV}(m)|_{\bar P_kV^{\rm reg} \cup \bar P_kV|_{O \se D}}$.
The complement of the latter set is of codimension two in $\bar P_kV|_{O}$,
so $s$
 extends to a section of ${\cal O}_{\bar P_kV}(m)|_{O}$.

Conversely let $\tilde Q$ be an $m$-linear meromorphic function on
${\cal O}_{\bar P_kV}(-1)|_{O}$ corresponding to a meromorphic section
of ${\cal O}_{\bar P_kV}(m)|_{O}$. Then $Q:= \tilde Q \circ \tilde \alpha_k$
is a meromorphic function on $\bar J_kV^{\rm reg}|_O$. By the Riemann
Extension
Theorem, it extends to a meromorphic function on $\bar J_kV|_U$. It satisfies
equation (\ref{pf5}) on $\bar J_kV^{\rm reg}|_O$ since $\tilde Q$ corresponds
to a section of ${\cal O}_{\bar P_kV}(m)|_{O}$ and since the
fibers of $\tilde \alpha_k$ are invariant under the action of $G_k^o$.
Hence, it satisfies
equation (\ref{pf5}) over $O$ by the Identity Theorem.

Finally we remark that if we start with a holomorphic rather than a meromorphic
function (respectively, section) in the arguments above, we would obtain a
holomorphic section
(respectively, function) as a result.
\qed

\section{Log-directed jet metrics}
\setcounter{equation}{0}
\subsection{The case of 1-jets}
This case was already treated in the second named author's thesis. We
recall the basic results
after some definitions.

For a line bundle $ L$ over a complex variety $ X$,
let $E_L$ be the union of the base locus
$$\mbox{Bs} |L| := \{x\in X:\ s(x)=0 \mbox{ for all }
s\in H^0(X, L)\},$$
of $L$ and
the restricted exceptional locus
$$\{x\in X\setminus \mbox{Bs} |L|\
:\ \dim_{x} \ph_L^{-1}(\ph_L(x)) >0 \}$$
of the rational map
$$\ph_L:= [s_1:...: s_n]: X \ \cdots \r \pr^{n-1},$$
where $\{s_1,...,s_n \}$ is a basis of $H^0(X, L)$.
  We will call $E_L$ the
{\em basic locus} of $L$.
 Define the {\em stable basic locus} of $L$ to be
$$
S _{L} : = \bigcap _{m>0} E _{mL} .
$$
A standard argument (worked out in details
in the appendix) shows that for any line
bundle $H$ on a normal variety $X$, we have
$$\mbox{Bs} |mL-H| \s S_L $$
for some sufficiently large $m$.
Let $X_{\rm reg}$ be the smooth part of $X$. Let $L^*$ be the dual bundle
of $L$. Recall that a continuous function
$g: L^* \r [ 0, \infty ]$ such that
\begin{equation} \label{m1}
g(cv) = |c|^2 g(v)
\end{equation}
for all $c \in \cz$ and $v \in L^*$ is
called a {\em singular metric} on $L^*$.
By equation (\ref{m1}), the set $g^{-1}(0)
\cup g^{-1}(\infty)$ consists of the zero section of $L^*$ and the inverse
image
in $L^*$ of a closed subset $\Sigma_g$ of $X$.
For our purpose, we will always assume that
the open set $U=X_{\rm reg} \se \Sigma_g$
is dense in $X$ and that $g$ is twice
differentiable on $L^*|_U$.
Then $dd^c \log g$ is a real $(1,1)$ form
outside the zero section of $L^*|_U$
 invariant under the $C^*$ action
given by equation (\ref{m1}) and is thus the
pull back of a real $(1,1)$ form on $U$
denoted by
$${\rm Ric} (g) = \Theta_{g^{-1}}
= \Theta_{g^{-1}}(L)\, ,$$
which is known as the curvature form of $g$.
By convention, $g$ is called a
{\em pseudometric} if $g^{-1}(\infty) =
\emptyset$ and $g$ is called a {\em metric}
if $\Sigma_g = \emptyset$.

Let now $(\bar X, D)$ be a log-manifold, and set again
$X = \bar X \setminus D$.
A K\"ahler metric $\omega$ on $X$ gives a metric on $TX$ which in turn gives a
metric $g_{\omega}$ on ${\cal O}_{\bar {P}_1 X}(-1)|_{P_1 X}$.
If $\omega$ behaves logarithmically
along $D$, then $g_{\omega}$ extends to a metric on
${\cal O}_{\bar {P}_1 X}(-1)$ which
we can use to dominate a scalar multiple of any pseudometics on
${\cal O}_{\bar {P}_1 X}(-1)$
 by appealing to the compactness of $\bar X$. This is the
basic strategy  used to obtain the following result
(Proposition 1 of \cite{Lu}). We remark that Noguchi (\cite{No0}) had already
similar results in the case $X$ is compact under the  assumption that
${\cal O}_{ P_1 X}(m)$ is
spanned by global sections everywhere
on ${\cal O}_{ P_1 X}(m)$ for $m$ large
enough.

\begin{prop}\label{Ric} Assume $(\bar X, D)$
is a log-manifold and that $\bar X$ is
projective.
Let $\pi_1: \bar P_1 X \r \bar X$ be the natural projection, let $\bar \Xi$
be a subvariety of
$\bar P_1 X$ and let
$\sigma:\bar Z \r \bar\Xi\seq \bar P_1 X$
be the normalization of $\bar \Xi$. Let
$\bar L_\sigma=\sigma^{-1}{\cal O}_{\bar P_1 X}(1)$,
$Z=\sigma^{-1}( P_1 X)$ and
$L_\sigma=\bar L_\sigma|_Z$.  Then there is a pseudometric
$g$ on $L_\sigma^*= \sigma^{-1}({\cal O}_{ P_1 X}(-1))$ with $\Sigma_g \s
S_{{\bar L_\sigma}}
$,
such that ${\rm Ric} (g)$ is the pullback of a K\"ahler metric $\omega$ on $X$,
specifically
$${\rm Ric} (g)=(\sigma \circ \pi_1)^* \omega\, ,$$
such that this K\"ahler metric
$\omega$ dominates  $g$, in the sense that
\begin{equation} \label{m2}
(\sigma^* g_{\omega})(\xi ) \geq
g (\xi )
\end{equation}
for all $\xi \in L_{\sigma}^*$ outside
$S_{\bar L_{\sigma}}$ and $\sigma^{-1}({\rm Sing}(\bar \Xi))$.
\end{prop}
By the usual definition of holomorphic
sectional curvature, we see that equation
(\ref{m2}) says precisely that $g$,
as a ``length'' function on $X$ in the tangent directions
defined by $\bar\Xi$,
has holomorphic sectional curvature
bounded from above by $-1$, and $g$ is nonvanishing outside
$S_{{\bar L_\sigma}}$.

Hence, the usual Ahlfors' Lemma applies to show that
if $f: \Delta \r X$ is any holomorphic map
 from the unit disk $\Delta \s \cz$ whose lifting $f_{[1]}$ has image in
$\bar\Xi$ but not completely in
$\sigma (S_{\bar L_\sigma}) \cup
{\rm Sing}(\bar \Xi)$, then
$f$  must satisfy the distance decreasing property.

{From} this, the following result is derived by elementary arguments in
\cite{Lu}
(see also Noguchi (\cite{No0})), which we quote.

\begin{theo}\label{main} With the same setup as in Proposition \ref{Ric},
we let $\Delta^* = \Delta \se \{0\}$ be the punctured unit disk and set
$\bar L_0= {\cal O}_{\bar P_1 X}(1)|_{\bar\Xi}$. \\
{\rm (a)} (Distance decreasing property)
If $f:\Delta \r X$ is a holomorphic map whose lift $f_{[1]}$ has
values in $\Xi$ but not entirely in $S_{\bar L_0}\cup {\rm Sing}(\bar
\Xi)$, then $f^* g \leq \rho$
(where $\rho$ is the Poincar\'{e} metric  on $\Delta$).\\
{\rm (b)} (Degeneracy of Holomorphic Curve)
If $f:\cz \r X$ is holomorphic such that $f_{[1]}$ has values in $\bar
\Xi$, then
$f_{[1]}(\cz) \s S_{{\bar L_0}}\cup {\rm Sing}(\bar \Xi)$.\\
{\rm (c)} (Big Picard Theorem) If $f:\Delta^* \r X$
is holomorphic such that $f_{[1]}$ has values in $\bar \Xi$ but not entirely in
$S_{{\bar L_0}}\cup {\rm Sing}(\bar \Xi)$, then $f$ extends to a
holomorphic map $\bar{f}: \Delta \r \bar{X}$.
\end{theo}

\subsection{The general case}
We call a pseudometric $h$ on
${\cal O}_{ {P}_k V}(-1)$
a {\it $k$-jet pseudometric} on $(X,V)$,
and define $B_k =S_{{\cal O}_{\bar {P}_k V}(1) }$, the stable basic locus of
${\cal O}_{\bar {P}_k V}(1)$.
\begin{theo} \label{ahlfors} With the notations as in Theorem \ref{main},
assume that $(\bar X,D,\bar V)$ is
a log-directed
manifold and that $\bar X$ is projective.
\begin{itemize}
\item[a)] If $B_k \not= \bar P_k V$, then
there exists a k-jet pseudometric $h$ on
$(X,V)$ with $\Sigma_h \seq B_k$
such that $h$ has curvature bounded from above by $-1$
in the sense that ${\rm Ric}(h)=\pi_k^* \omega$ is the pullback of a
K\"ahler metric $\omega$ on
${P}_{k-1} V$ such that $g_{\omega}$ dominates $h$. In particular, we have
$$< \Theta_{h^{-1}},|\xi|^2>\ =\ <{\rm Ric}(h),|\xi|^2> \ \,
\geq \, h((\pi_k)_* \xi
) \ \ \ {\rm for}\ \ \xi \in V_k .$$
\item[b)] If $f: \cz \r \bar X \setminus D$
is holomorphic with
$f_*(T \cz) \s \bar V$, then $f_{[k]}(\cz ) \s B_k$.
\item[c)]  If $f: \Delta^* \r \bar X \setminus D$ is holomorphic with
$f_*(T \Delta^*) \s \bar V$, then:\\ Either $f$ extends to a
holomorphic map $\bar f: \Delta \r \bar X$ or  $f_{[k]}(\Delta^* ) \s B_k$.
\end{itemize}
Moreover, let $Y \s \bar P_k V$ be any subvariety.
We define $B_k(Y)=S_{{\cal O}_{\bar {P}_k V}(1)|_Y}$.
If $f$ lifts to a map with values in $Y$, then
b) and c)
hold with $B_k(Y) \cup {\rm Sing}(Y)$ instead of $B_k$.
\end{theo}
\noindent {\bf Proof }
Apply Proposition \ref{Ric} and Theorem \ref{main} to the log-manifold
$(\bar P_{k-1}V,$ $ D_{k-1})$
and the subvariety $\bar \Xi = \bar P_kV$ (or
$\bar \Xi = Y \subset \bar P_kV$).
  Note that,
since  $\bar V_k \s \bar T(P_{k-1}V)$ is a
 holomorphic subbundle,
 $\bar{ P}_1 V_k = \bar {P}_k V$ is a submanifold in $\bar{ P}_1(P_{k-1}V)$
and ${\cal O}_{ {P}_k V}(-1)=
{\cal O}_{\bar{ P}_1({P}_{k-1} V)}(-1)|_{ {P}_k V}$ by Proposition \ref{subm}.
 \qed

\section{Logarithmic Bloch's  and Lang's Conjecture}

\setcounter{equation}{0}

In this section we apply our method to the special case of
semi-abelian varieties where our Ahlfors-Schwarz Lemma
(Theorem~\ref{ahlfors}) gives
a logarithmic version of Bloch's Theorem and our big Picard Theorem
yields a big Picard version of Bloch's Theorem. By using the
Wronskian  associated to the theta function of an effective
divisor in a semi-abelian variety
(\cite{SY}, \cite{PNog}), we affirm furthermore a logarithmic
version of Lang's Conjecture and a big Picard analogue of it,
all via metric geometry on Logarithmic Demailly-Semple jets.

\subsection{Statement of the results}

We first recall the definition and some basic facts on
semi-abelian varieties (see \cite{PNog}, \cite{Ii1}, \cite{Ii2})
needed to state our results.

A quasiprojective variety $G$ is called a semi-abelian variety
if it is a commutative group which admits an
exact sequence of groups
$$0 \r (\cz^*)^{\ell} \r G \r A \r 0,$$
where $A$ is an abelian variety of dimension $\sf m$.

 Taking the pushforward of $(\cz^*)^{\ell}\s G$ with the natural embedding
$(\cz^*)^{\ell} \s (\pr^1)^\ell$, we obtain a smooth completion
$$\bar{G} = (\pr^1)^{\ell} \times_{(\cz^*)^{\ell}}  G$$
of $G$ with boundary divisor $S$, which has only
normal crossing singularities. We denote the natual action of $G$ on
$\bar{G}$ on the right as addition.
It follows that the exponential map from the Lie algebra $\cz^n$
is a group homomorphism and, hence, it is also the universal covering map of
$G=\cz^n / \Lambda$, where $\Lambda=\Pi_1(G)$
is a discrete subgroup of $\cz^n$ and $n={\sf m}+\ell$.

Following Iitaka (\cite{Ii2}), we have the following trivialization of
the logarithmic tangent respectively cotangent bundles of $\bar{G}$:
Let $z_1, ...,z_n$ be the
standard coordinates of $\cz^n$. Since $dz_1,...,dz_n$ are invariant under
the group action of translation on $\cz^n$, they descend to
forms on $G$. There they extend to
logarithmic forms on $\bar G$ along $S$, which are elements of
$H^0(\bar G, \bar \Omega G)$. These logarithmic $1$-forms
are everywhere
linearly independent on $\bar G$.
Thus, they
globally trivialize the vector bundle
$\bar\Omega G$.
Finally, we note that these logarithmic forms
are invariant under the group action of $G$ on $\bar{G}$, and, hence,
 the associated trivialization of $\bar\Omega G$ over $\bar G$
is also invariant.

We now state the main theorems of this section. With the above
setup, let  $f:\Gamma\rightarrow G$  be a holomorphic map,
where $\Gamma$ is either $\cz$ or the punctured disk $\Delta^*$.
Denote by $\bar X(f)$ the Zariski closure of $f(\Gamma)$ in $\bar G$ and let
$X(f)=\bar X(f)\cap G$.
Furthermore, let
$D \s G$ be a reduced algebraic divisor in $G$,
which we regard as a union of codimension
one algebraic subvarieties of $G$. We note that an algebraic subvariety of
$G$ which is also a subgroup is necessarily a semi-abelian variety as well,
see \cite{Ii2}.

\begin{theo} \label{c}$ $\\[-.3in]
\begin{itemize}
\item[(a)] Let $f: \cz \r  G$ define a holomorphic curve.
 Then $X(f)$ is
a translate of an algebraic subgroup of $G$.

\item[(b)] Let $f: \cz \r (G \setminus D)$ be holomorphic.
Then $X(f) \cap D = \emptyset$.
\end{itemize}
\end{theo}

\begin{cor} \label{cc}
If $D$ has nonempty intersection with any translate of an algebraic
subgroup of $G$ (of positive dimension), then $G \setminus D$ is
Brody hyperbolic.\\
In particular, this holds if $G=A$ is an abelian
variety and $D$ is ample. \qed
\end{cor}

\noindent  Theorem \ref{c} (a) is a logarithmic version of
Bloch's Theorem, first proved by Noguchi (\cite{No2}),
 (b) is a logarithmic version of Lang's Conjecture.
Both Theorem \ref{c}  and Corollary \ref{cc}  were obtained by
Noguchi (\cite{PNog}), and were, in the nonlogarithmic case,
first proved by  Siu-Yeung (\cite{SY}).

\begin{theo} \label{bp}$ $\\[-.3in]
\begin{itemize}
\item[(a)] Let $f:\Delta^* \r G$ be a holomorphic map. Then either it extends
to a holomorphic map $\bar{f}:\Delta \r \bar{G}$, or there exists a
maximal algebraic subgroup
$G'$ of $G$ of positive dimension such that  $X(f)$ is foliated by
translates of $G'$.

\item[(b)] Let $f: \Delta^* \r (G \setminus D)$ be holomorphic.
Then one of the following holds:\\
(i) $f$ extends to $\bar{f}: \Delta \r \bar{G}$.\\
(ii) $X(f) \cap D = \emptyset$.\\
(iii) There is an algebraic subgroup $G'' \s G'$ of positive
dimension such that
$X(f) \cap D$  foliated by translates of $G''$.

\item[(c)] Let now $f: \Delta^* \r (A\setminus D)$, where $G$ is the
special case of an abelian variety $A$.
Then one of the following  holds:\\
(i) $f$ extends to $\bar{f}: \Delta \r A$.\\
(ii) There exists an algebraic subgroup $G'' \s G'$ of positive dimension
such that $D$ is foliated by translates of $G''$.
\end{itemize}
\end{theo}

\begin{cor} \label{cbp}
If $G=A$ is an abelian variety and $D$ is ample, then
$f:\Delta^* \r A\setminus D$
extends to a holomorphic map $\bar{f}:\Delta \r A$.
\end{cor}

\noindent We remark that Theorem \ref{bp} and
 Corollary \ref{cbp} are big Picard type Theorems.
Aside from  (a), which can be found in Noguchi (\cite{No2}),
these are, to our knowledge, new to the literature. \\

\noindent {\bf Proof of Corollary \ref{cbp} }
Corollary \ref{cbp} follows from  Theorem \ref{bp} (c) and the fact that an
ample divisor $D$ in an abelian variety $A$ cannot be foliated
by translates of an algebraic subgroup  $A''$ of $A$ of positive dimension.
For assume it were.
Then $D= q^{-1}(\bar D)$, where $\bar D$
is a divisor in $A/A''$ and  $q:A \r A/A''$
is the quotient map. But then ${\cal O}(D)=q^{-1}{\cal O}(\bar D)$ is
trivial along $A''$, since $A''$ is
a fiber of the map $q$. This is a contradiction.
\qed $\ $

\noindent {\bf Remark } The last part of Corollary \ref{cc} follows from
Corollary \ref{cbp} as follows.
The $\cz\setminus \Delta$ is biholomorphic to $\Delta^*$.
So we can conclude
from Corollary \ref{cbp} that any entire curve $f:\cz \r A$ extends
to a holomorphic map $\bar{f}:\pr^1 \r A$. Hence,  $\bar f$  must be constant,
since all coordinate 1-forms on $A$ must pull back to the zero on $\pr^1$ as
$\pr^1 $ has no nontrivial 1-forms.

 However, Corollary \ref{cbp} does not follow from Corollary \ref{cc}. It
would if  $A \se D$ were hyperbolically embedded in $A$. This
would be the case, for example,
if $D$ were hyperbolic (see for example \cite{La1}).
But even a very ample divisor in $A$ is not hyperbolic in general.
To see this,
choose any translate $T$ of an algebraic subvariety which is of codimension
at least 3 in $A$.
Then there always exists an irreducible ample divisor
in $A$ which contains $T$, as can be deduced by applying the
Bertini's Theorem 7.19 in \cite{Ii}.
Note that a hyperbolic open subset $V$ in a projective variety
$\bar V$ need not be hyperbolically embedded in general, as one can
easily see by blowing up a point of $\bar{V} \setminus V$.\\

\noindent {\bf Remark } Let $G = (\cz^*)^n \s \pr^n$, $n \geq 4$.
Let $\bar{D}$ be a generic
hyperplane in $\pr^n$ and $\bar{H}$ be another hyperplane with
 $G \cap \bar{H} \not= \emptyset$ and
$\bar{H} \cap \bar{D} \cap G = \emptyset$. Then $\bar{H} \cap (G \se \bar{D}
)$ is equal to $\pr^{n-1}$ minus at most $n+2$ hyperplanes, which contains
nontrivial images of $\cz$ and hence, admits maps $f$ from $\Delta^*$
which do not extend to $\Delta$. This is because the complement of
$n+2$ hyperplanes in $\pr^{n-1}$ contains nontrivial diagonals for $n \geq 4$,
which are nonhyperbolic. So we get
examples of $f$
for Theorems \ref{c} (b) and  \ref{bp} (b) (ii) with
nontrivial $X(f)$.\\

\noindent {\bf Remark } Let $A$ be an abelian variety,
$D \s A$ a divisor and $f_1:\Delta \r A\setminus D$  a holomorphic map.
Let $X(f_1)$ be the Zariski closure of $f_1(\Delta) \s A$. Let $E$ be an
elliptic curve and $q:\cz \r E$ be the universal cover. Then
$$f(z)=(f_1(z), q \circ {\rm  exp}(\frac{1}{z})): \Delta^* \r A\times E$$
does not extend.\\

\noindent This easy construction provides examples which are relevant to
Theorem \ref{bp}:
\begin{itemize}
\item[(1)] It makes  Theorem \ref{bp} (b) (iii) and  (c) (ii) sharp.
\item[(2)] Choose $f_1$ in such a way that $X(f_1)$ is not a translate of an
algebraic subgroup in $A$. Then we have an example for
 (a) where $X(f)$ is not itself a translate of an algebraic subgroup
of $A\times E$.
\end{itemize}

\subsection{Some results on semi-abelian varieties}
We first summarize some elementary properties of semi-abelian varieties.

\begin{lem} \label{sav}$ $\\[-.3in]
\begin{itemize}
\item[(a)] The quotient of a semi-abelian variety $G$ by an algebraic subgroup
$G'$ is again a semi-abelian variety,
and the quotient map $q:G \r G / G'$ is an
algebraic morphism.
\item[(b)] If $X \s G$ is an algebraic variety foliated by translates of
$G'$,
then $$X/G' \s G/G'$$ is again an algebraic variety.
\item[(c)] If $X$ is an algebraic subvariety of $G$, and $h:X \r \pr^1$
is a rational function,
then the closed subgroup
$$G'= \{ a \in G : X=(X+a)\} \cap \{ a \in G : h(x) =h(x+a) \;
{\rm for\ all}\; x \in X \}$$
is again an algebraic subvariety.
\end{itemize}
\end{lem}

\noindent {\bf Proof }
 Lemma \ref{sav} should be well known, but since we do not know
a precise reference we indicate a proof. {From} the fact that connected
algebraic subgroups of a semi-abelian variety are again semi-abelian,
it is easy to see that  one
can consider quotients of $G$ by
$G'$ by taking the quotient on the abelian and the $({\cz^*})^{\ell}$ factors
separately. Now the quotient of the abelian factor by an algebraic subgroup
is abelian by isogeny, and the quotient of $({\cz^*})^{\ell}$ by a
connected algebraic subgroup is likewise a product of $\cz^*$,
see \cite{Ii2}. Hence, we obtain a $({\cz^*})^l$ bundle over
an abelian variety for some $l$, which projectivizes to a $\pr^l$
bundle over
a projective variety, and, therefore, must be projective. {From} this the
entire Lemma \ref{sav} follows. \qed

\begin{lem} \label{proj}
Let $A$ be an abelian variety and $D \s A$ a  reduced algebraic divisor.
Let $A'$ be an algebraic subgroup of $A$ and $T$ a translate of $A'$
in $A$.
Assume $T \cap D = \emptyset$. Then $D$ is foliated by translates of $A'$.
\end{lem}

\noindent {\bf Proof }
Without loss of generality we may assume that $D$ is irreducible.
Let $q: A \r A/A'$ denote the quotient map. Since $q$ is a proper map,
$q(D)$ is a projective subvariety in $A/A'$. Since $D$ is irreducible
and $T \cap D = \emptyset$, $q(D)$ is an irreducible divisor. So
$\tilde{D} = q^{-1}(q(D)) \s A$ is also an irreducible
divisor containing $D$ as $q$ is smooth.
This forces  $D=\tilde{D}$.\qed
$\ $

\noindent {\bf Remark } Lemma \ref{proj} is  false for
semi-abelian varieties. For we may take
$G=(\cz^*)^2$, $T=\{(z_1,z_2) \in G: z_1=1 \}$,
$D=\{(z_1,z_2) \in G: (z_1-1)z_2 =1 \}.$

\subsection{Jet bundles on semi-abelian varieties}

To simplify notation,
we work exclusively with a semi-abelian variety $G$
and its associated log-manifold $(\bar{G},S)$ defined as before.
We remark, however, that many the definitions and results
hold for an arbitray log-manifold.

Recall from  subsections 1.2 and 2.1 that $\bar{P}_kV$
denotes the  logarithmic $k$-jet bundle
of $(\bar{G}, S, \bar V)$ and that $\bar P_kG$
denotes the logarithmic jet bundle of $(\bar{G}, S)=(\bar{G}, S, \bar{T}G)$.
Note that the
log-directed morphism $i:(\bar{G}, S, \bar V) \r (\bar{G}, S, \bar{T}G)$
induces a canonical realization of $\bar{P}_kV$ as a submanifold of
$\bar P_kG$ as $\bar V$ is a subbundle of $\bar TG$ over $\bar G$.

Let $D \s G$ be a
reduced algebraic divisor. Then, by Hironaka (\cite{H}),
there exist a log-manifold $(\bar{Y},E)$ and  a log-morphism
$p:(\bar{Y},E) \r (\bar{G}, S)$ with
\begin{itemize}
\item[1)] $p^{-1}(S \cup D) = E,$

\item[2)] $p: p^{-1}(\bar{G} \se D) \r
 \bar{G} \se D$
is biholomorphic.
\end{itemize}

\noindent
Given a subbundle $\bar V$ of $\bar TG$,
we have the following commutative diagram:
\begin{equation}
 \label{plp}
\begin{array}{ccccc}
\co_{\bar P_k ~Y} ~(-1)
~~ & ~ \stackrel{ ( p_{[k - 1]} )_{*} }{ \lra }
~~ & \co_{\bar P_k ~G} (-1)
~~ & \stackrel{ (i_{ [k - 1] })_{*} }{ \hookleftarrow }
~~ & \co_{ \bpk ~(V) } ~(-1) \\
& & & & \\
\da
&
& \da
&
& \da \\
& & & & \\
\bar P_k ~Y
& \stackrel{ p_{[k]} }{ \lra }
& \bar P_k ~G
& \stackrel{ i_{[k]} }{ \hookleftarrow  }
& \bpk ~(V) \\
& & & & \\
\da
&
& \da
&
& \da \\
& & & & \\
(~ \ol{Y}, E ~)
& \stackrel{ p }{ \lra }
& (~ \ol{G}, S ~)
& \stackrel{ i }{\hookleftarrow }
& (~ \ol{G}, S, \bar V ~)
 \end{array}
\end{equation}

\noindent Here, $i_{[k]}$ realizes
$\bar{P}_kV$ as a submanifold of $\bar P_kG$ by
Proposition~\ref{functoriality}.
Outside $\pi_{0,k}^{-1}(D) \s \bar P_kG$, the maps
$p_{[k]}$ (and hence, ${p_{[k-1]}}_*$) are isomorphisms.
All other maps are holomorphic.
We define
$$  \bar Z_k:=\overline{p_{[k]}^{-1}(\bar{P}_kV \se
\pi_{0,k}^{-1}(D))}^{\rm Zariski} \s \bar P_kY.$$

\begin{defi} \label{deflp}
A meromorphic section $s$ of ${\cal O}_{\bar{P}_kV}(m)$
is said to have at most log-poles along $D$
if it pulls back,
via the map $(p_{[k-1]})_*|_{\bar Z_k}=(p_{[k-1]}|_{\bar Z_k})_*$,
to a holomorphic
section\footnote{With this we mean
that, after pulling back the section $s$ over
the part
of $ \bar Z_k$ where the meromorphic map
$(p_{[k-1]})_*$ is holomorphic, it extends to a holomorphic
section of ${\cal O}_{\bar{P}_kY}(m)|_{\bar Z_k}$.}  of
${\cal O}_{\bar{P}_kY}(m)|_{\bar Z_k}$.
\end{defi}

\noindent
Suppose that a meromorphic section $s$ of ${\cal O}_{\bar{P}_kV}(m)$
over an open subset
$U\seq \bar G$ is defined  by a meromorphic function $Q$ on
$J_k G|_U$ satisfying equation~(\ref{invj}),
see Proposition \ref{dlog}. Suppose more precisely that
$Q$ is given by a polynomial
in the differentials up to order $k-1$ of sections of
$\bar \Omega G|_U$ as well
as the differentials up to order $k-1$ of $d\log \theta$, where
$\theta$ is a meromorphic function on $U$
nonvanishing and holomorphic on $U_0=U\se \{D\cup S\}$.
Then, after composing with $p$ given
above, these differentials, and so also  the
polynomial in them, become holomorphic functions
on $\bar  J_k Y|_{p^{-1}(U)}$ by Proposition~\ref{p1}(c).
Furthermore, the resulting polynomial still satisfies
equation (\ref{invj}) on $p^{-1}(U)$.
Hence, by Proposition \ref{dlog},
we obtain a holomorphic section of ${\cal O}_{\bar{P}_kY}(m)|_{\~U}$
that matches with the pullback of $s$ on an open set of $\bar Z_k$.
Therefore, $s|_U$ is meromorphic with at most log-poles along $D$. If
such a description is possible on a neighborhood $U$ of each point in $\bar G$,
then $s$ is meromorphic with at most log-poles along $D$. We will consider
examples of such an $s$ in subsection 5.5.

\begin{lem} \label{jj}
Let $(\bar{G},S,\bar V)$ be as above.
\begin{itemize}
\item[(1)] There exist  injective maps
$$\bar{P}_{k+l}V \r \bar P_l({P}_kV)\ \ \ \mbox{and }\ \
{\cal O}_{\bar{P}_{k+l}V}(-1) \r
{\cal O}_{\bar P_l(P_kV)}(-1)$$
which are given outside $S$
by $f_{[k+l]} \mapsto (f_{[k]})_{[l]} $ and by
$f_{[k+l-1]}' \mapsto (f_{[k]})_{[l-1]}'$, respectively, and
which realize $\bar{P}_{k+l}V \s \bar P_l({P}_kV)$ and
${\cal O}_{\bar{P}_{k+l}V}(-1) \s {\cal O}_{\bar P_l({P}_kV)}(-1)$,
respectively, as submanifolds.

\item[(2)] Furthermore, let $\Gamma$ be a curve and $f:\Gamma \r X$
be a holomorphic map
which is  tangent  to $V$.
As before, let $\bar X_k(f) \s \bar{P}_kV$ be the Zariski closure of the
image of
the $k$-th lift $f_{[k]}: \Gamma \r \bar{P}_kV$ of the map $f$.
We denote $\pi_{0,k}^{-1}(S)$
again by $S$ and  $\pi_{0,k}^{-1}(G) \cap \bar X_k(f)$ by $X_k(f)$. We
recall that by Hironaka there exists a log-morphism
$$\Psi: (\bar{\tilde X_k}(f), \tilde{S}) \r
(\bar{P}_kV, S)\ \ \ \ \mbox{such that:}$$
\begin{itemize}
\item[(a)] $\Psi (\bar{\tilde{X_k}}(f)) =   \bar X_k(f)$.
\item[(b)] $\Psi^{-1}(S) = \tilde{S}$.
\item[(c)] $\Psi$ is biholomorphic outside
$\Psi^{-1}\Big({\rm Sing} (\bar X_k(f))\Big)$.
\end{itemize}
We set $\tilde{X}_k(f) = \bar{\tilde{X}_k}(f) \se \tilde{S}$.
With this setup, we have the following commutative diagram:
\begin{equation}
\begin{array}{ccccc}
\bpl ~(~  {\tilde{X}_k} (f) ~)
~~ & ~ \stackrel{ \Psi_{[l]} }{ \lra }
~~ & \bpl ~(~ P_k ~V ~)  \\
& & \\
|
&
&  i \ua \\
& & \\
|
&  \bar X_{k + l} ~(f) ~ \subset
&  \bpkl ~(V) \\
& & \\
\da
& \da
& \da \\
& & \\
\bar{\tilde X_k}(f) ~~~~ \stackrel{ \Psi }{ \lra }
& \bar X_{k} ~(f) ~~~~ \subset
& \bpk  ~(V)
 \end{array}
\end{equation}

where $\Psi_{[l]}$ may be meromorphic, all other maps in
the diagram are holomorphic and the following holds:
$$i(\bar X_{k+l}(f)) \s \Psi_{[l]} (\bar P_l(\tilde{X}_k(f))).$$

\item[(3)] Let $s$, $t$ be meromorphic sections of the bundle
${\cal O}_{\bar{P}_kV}(m)$
with at most log-poles along $D$, and assume $t$ is not the zero section.
 Then $t^2 \cdot d(\frac{s}{t})$ can be considered as a
meromorphic section of the bundle
${\cal O}_{\bar{P}_{k+1}V}(2m+1)$ with at most log-poles along $D$.
\end{itemize}
\end{lem}

\noindent  {\bf Proof for (1) }
This follows directly from Corollary~\ref{subm} and
Proposition~\ref{functoriality} applied to the subbundle
inclusion $\bar V_k \s \bar T P_{k-1} V$. The presentation of
the maps outside $S$ follows from Proposition~\ref{Fof}.

\noindent {\bf Proof for (2) }
$\Psi_{[l]} : \bar P_l(\tilde{X}_k(f)) \r \bar P_l ( {P}_kV)$
is a proper rational
map. Hence,
$\Psi_{[l]}(\bar P_l(\tilde{X}_k(f))$ is an algebraic subset
containing $i(f_{[k+l]}(\Gamma)) = (f_{[k]})_{[l]}(\Gamma)$.
Therefore, we also have
$i(\bar X_{k+l}(f))$ $ \s \Psi_{[l]} (\bar P_l(\tilde{X_k}(f))$.

\noindent {\bf Proof for (3) }
 $\frac{s}{t}$ is a rational function on $\bar{P}_kV$.
Hence, $d(\frac{s}{t})$ is a rational section of
${\cal O}_{\bar P_1({P}_kV)}(1)$,
which lifts back, via the inclusion map
${\cal O}_{\bar{P}_{k+1}V}(-1) $ $\r
{\cal O}_{\bar P_1({P}_kV)}(-1)$, to a rational section of
${\cal O}_{\bar{P}_{k+1}V}(1)$.
Hence, $t^2 \cdot d(\frac{s}{t})$ is a rational section of
${\cal O}_{\bar{P}_{k+1}V}(2m+1)$.
To prove that $t^2 \cdot d(\frac{s}{t})$ again is a meromorphic
section with at most log-poles
along $D$, we pull back the sections $s$, $t$  and $t^2 \cdot d(\frac{s}{t})$
to $\bar Z_k$ by the map $(p_{[k-1]})_*$.
Then by definition the sections $s$ and $t$ become holomorphic sections
on $\bar Z_k$. It suffices to show that the rational section
$t^2 \cdot d(\frac{s}{t})$ also is holomorphic on $\bar Z_k$.
It suffices to prove this locally on $\bar Z_k$. Given any point $w \in
\bar Z_k
\s \bar P_kY$ there exists an open neighborhood $U$ of $w$ in $\bar P_kY$
such that the holomorphic sections $s$ and $t$ over $\bar Z_k \cap U$
extend to holomorphic sections of the bundle ${\cal O}_{\bar P_kY}(m)$ over
$U$,
and that this bundle is trivial over $U$. After choosing such a trivialization,
one has, by the product rule for the holomorphic functions $s$ and $t$,
$$t^2 \cdot d(\frac{s}{t}) = tds-sdt\, ,$$
the latter being a holomorphic section of ${\cal O}_{\bar P_1(U \se E)}(1)$.
 \qed

\vspace{0.3cm}

Next, we want to prove that $\bar P_kG$, and even $\bar{P}_kV$,
are trivial over $\bar{G}$
for certain subbundles $\bar V \s \bar{T}G$, which we will call special.
 Let
$z_1$, ..., $z_n$ be any linear coordinates of the universal cover $\cz^n
\r G$.
 We have observed that $\bar{T}G = \cz^n \times \bar{G}$,
where the trivialization is given
by $dz_1$, ... $dz_n$.
\begin{defi}
$\bar V$ is said to be special if $\bar V = \cz^r \times \bar{G}$
in this trivialization.
\end{defi}
{From} now on, all subbundles $\bar V \s \bar{T}G$ we use
are assumed to be special.

\begin{lem} \label{sat} $ $\\[-.3in]
\begin{itemize}

\item[(a)]  The map
\begin{equation} \label{iso} {\cal O}_{P_kV}(-1) \r ({\cal O}_{P_kV}(-1))_0
\times G;\:\:
f_{[k-1]}'(0) \mapsto ((f-f(0))_{[k-1]}'(0), f(0))
\end{equation}
gives an isomorphism  ${\cal O}_{P_kV}(-1) \r {\cal
O}_{\bar{P}_kV}(-1)|_G$, and
this isomorphism is invariant under the action of $G$.
\item[(b)] This isomorphism  can be extended to a $G$-invariant isomorphism
$$ \bar{\Psi}_{k} : {\cal O}_{\bar{P}_kV}(-1) \r ({\cal O}_{P_kV}(-1))_0
\times \bar{G}$$
respecting the fibers of the line bundles.
\item[(c)] By combining with the canonical line bundle projections we get the
same isomorphisms
with $P_kV$ and $\bar P_kV$ instead of $ {\cal O}_{P_kV}(-1)$ and
${\cal O}_{\bar{P}_kV}(-1)$.

\end{itemize}
\end{lem}

\noindent {\bf Proof }
 (a) is immediate, and (c) follows immediately from (b).
To prove  (b),
we use that by Corollary \ref{coor} (a) the trivialization of (a) extends
locally, so by the Identity Theorem it extends globally. The invariance under
the action of $G$ of this trivialization extends from $G$ to $\bar G$
by continuity.

We would like, however, also to indicate a direct proof.
It is obtained by proving the following more precise statement
by induction over $k$.\\

\noindent {\bf Claim S(k)}
Then there exist  trivialization maps $\bar{\psi}_{k}$,
induced canonically by the trivialization of $V$ by $dz_1,...,dz_n$, such that

$$
\begin{array}{cccccc}
& \bar V_k
~~ & ~ \stackrel{ \ol{\psi}_{k} }{ ~ \lra }
~~ & ( V_k)_{0} \times \ol{G}  \\
& & & \\
(+)_k & \da  \pi_{k}
&
&  ~~ \da ( \pi_{k} )_{0} \times  id_{ \ol{G} } \\
& & & \\
& \bar{P}_{k-1} ~V
&  \stackrel{ \ol{\Psi}_{k - 1} }{ \lra }
&  (P_{k - 1} ~V)_{0} \times \ol{G}

 \end{array}
$$
commutes and the upper line projectivizes to
$$
\begin{array}{cccccc}
& \bar P_kV
~~ & ~ \stackrel{ \ol{\Psi}_{k} }{ ~ \lra }
~~ & ( P_kV)_{0} \times \ol{G}  \\
& & & \\
(++)_k & \da  \pi_{k}
&
&  ~~ \da ( \pi_{k} )_{0} \times  id_{ \ol{G} } \\
& & & \\
& \bar{P}_{k-1} ~V
&  \stackrel{ \ol{\Psi}_{k - 1} }{ \lra }
&  (P_{k - 1} ~V)_{0} \times \ol{G}

 \end{array}
$$
where $\Psi_k$ extends the isomorphism in equation (\ref{iso}) in a
$G$-invariant
way.\\

\noindent Now {\bf S(1)} is clear from the trivialization
$$ V\ \ \tilde{\lra}\ \ \cz^r \times \bar{G},$$
since we are given that $V$ is special.
Assuming, by induction, that {\bf S(k)} is true.
Then we get $(++)_{k+1}$ by projectivizing
$(+)_{k+1}$. It remains to show
$(+)_{k+1}$.
By $(++)_k$ we get induced trivializations of the logarithmic tangent
bundles
$$
\begin{array}{ccccc}
 \bar T (P_kV)
~~ & ~ \stackrel{ (\ol{\Psi}_{k})_ *}{ ~ \lra }
~~ & T((P_kV)_{0}) \times \ol{T}G & \r & T((P_kV)_{0}) \times \cz^r
\times G\\
 & & & &\\
 \da  (\pi_{k})_*
&
&  ~~ \da (( \pi_{k} )_{0})_* \times  (id_{ \ol{G} })_*  & &
\\
& & & &\\
 \bar T({P}_{k-1} ~V)
&  \stackrel{ (\ol{\Psi}_{k - 1})_* }{ \lra }
& T(( P_{k - 1} ~V)_{0}) \times \ol{T}G & \r &
   T((P_{k-1}V)_{0}) \times \cz^r \times G\\

 \end{array}
$$
where the isomorphisms on the right hand side are obtained by trivializing
$\bar T G$ by the forms $dz_1,...,dz_n$.
We want to show it we restrict the isomorphism  in the upper line of this
diagram from $\bar T (P_kV)$ to $\bar V_k$, we also get a trivialization of
$\bar V_k$ over $\bar G$. Then we can denote this trivialization by
$\psi_{k+1}$
and the rest follows easily.
The key point of the proof is now that by equation (\ref{O1}), namely
$$\bar V_k({G}) := (\pi_{k})_*^{-1}({\cal O}_{\bar{P}_kV}(-1))
\s \bar{T}(P_kV),$$
 the subbundle
$\bar V_k \s \bar T P_kV$ is defined in an intrinsic way which is compatible
with the isomorphisms of the diagram above. \qed

\begin{lem} \label{id}
Let $Y$ be a complex manifold, $Z \s Y$ a complex submanifold and denote
by $i:(Z,TZ) \r (Y,TY)$ the directed inclusion map.
Let $g: \Delta \r Y$ be holomorphic with $dg(0) \not= 0$,
where $\Delta$ denotes again the unit disk in $\cz$.
Assume that $g_{[l]}(0)$ is in the image of the (composed) morphism
\begin{equation} \label{pf100}
 P_lZ \stackrel{ i_{[l]} }{ \lra } i^{-1}P_lY \r P_lY
\end{equation}
for all $l \geq 0$.
Then $g(\Delta) \s Z$.
\end{lem}

 \noindent {\bf Proof }
By Proposition \ref{functoriality} a) and b) the map in equation (\ref{pf100})
is a morphism.
 Let $U\s Y$ be a neighborhood of $g(0)$
 and $F:U \r \cz$ be a holomorphic function with $F|_{Z \cap U} \equiv 0$.
 It suffices to show $F \circ g \equiv 0$.
 There exist a small disk $\Delta_{\epsilon}$ and
 a map $h_l : \Delta_\epsilon \r Z$
 with $i_{[l]}((h_l)_{[l]}(0))=g_{[l]}(0)$ and by Corollary
\ref{corfunctoriality}
we may assume $dh_l(0) \not= 0$. By Proposition \ref{Fof} we have
$i_{[l]}(h_l)_{[l]}=(i \circ h_l)_{[l]}$ and hence,
$(i \circ h_l)_{[l]}(0)=g_{[l]}(0)$ and $d(i \circ h_l)(0) \not= 0$.
 Hence, we can reparametrize $i \circ h_l$ in a way that it has the same
 Taylor expansion as $g$ up to order $l$. We assume that this has been
 done, and note that $h_l$ still maps a neighborhood of the origin to $Z$.
 Hence,
 $$(\frac{\partial^j}{\partial t^j} F \circ g)(0)=
 (\frac{\partial^j}{\partial t^j} F \circ i \circ h_l)(0)
 =0 \;\;\;{\rm for} \;\;\; j \leq l.$$
Since $l$ is arbitrary, we get $F \circ g \equiv 0$. \qed

\subsection{The Main Lemma}

The following Main Lemma is the key step in proving Theorem \ref{c} and
Theorem
\ref{bp}. In the case $\Gamma = \cz$, it is a generalization of a lemma
contained in
\cite{PNog}, and for $G=A$,  it generalizes a lemma in \cite{SY}.

\begin{mlem} \label{ji}
With the same setup as that for Theorem \ref{c} (or  Theorem \ref{bp}),
let $\bar V \s \bar{T}G$ be a special subbundle. Assume that
$f: \Gamma \r G\se D$ is tangent to $V$ and, in the case $\Gamma = \Delta^*$,
that $f$ does not extend to $\Delta$ as a map to $\bar G$.
Let, for $k \geq 0$, $\bar X_k(f)$ denote
the Zariski closure of $f_{[k]}(\Gamma)$ in $\bar{P}_kV$.
Let, for $k, m \geq 1$, $\Theta$ be
meromorphic section of the line bundle ${\cal O}_{\bar P_kV}(m)$ with at most
log-poles along $D$, there exists an algebraic subgroup
$G' \s G$, of positive dimension, which leaves $\bar X_k(f)$ and, for
$k \geq 1$, also
$\Theta|_{\bar X_k(f)}$ invariant.

\end{mlem}

\noindent {\bf Remark }
The same is true for finitely many different sections $\Theta$.\\

\noindent
The rest of this subsection will be devoted to the proof of the Main Lemma.
It suffices to consider the case $k \geq 1$. In fact, to prove the case
$k=0$
we apply
the Main Lemma for $k=1$ and for
$\Theta$ being the zero section
in ${\cal O}_{\bar P_1V}(1)$.
Since the map $\pi_1: \bar P_1V
\r \bar X$ is equivariant under the action of $G'$ and maps $\bar X_1(f)$
surjectively $\bar X_0(f)$,
the subgroup $G'$ also leaves
$\bar X_0(f)$ invariant. So for the rest of the proof of the Main Lemma
we assume $k \geq  1$.\\

 We fix $u \in \Gamma$ to be any point for which $df(u) \not= 0$.
Then all $f_{[k+l]}(u)$, $l \geq 0$, are regular jets.
Let $s_0 \in H^0(\bar{P}_1V,{\cal O}_{\bar{P}_1V}(1))$ be a global section,
which is invariant under the action of $G$, and which is nonvanishing at
$f_{[1]}(u)$. It exists because $df(u) \not= 0$ and
${\cal O}_{\bar{P}_1V}(1) = {\cal O}_{P_1V_0}(1) \times \bar{G}$
(see Lemma \ref{sat}).
Choose an infinite sequence $\{n_0,\,n_1,\,n_2,\,n_3,...\}$
of natural numbers such that  the
following two conditions hold:
$$(2(k+l)-1) | n_l\; {\rm for}\; l \geq 0,$$
$$n_l \geq 2(n_{l-1}+1)\; {\rm for}\; l \geq 1.$$
For example, $n_l=2^l(2(k+l)-1)m$,
where $m={\rm deg} \Theta$ will work.
Let
 $$\Theta_0 = \Theta \cdot (s_0^{n_0 - m})$$
and, for $l \geq 1$, define inductively:
$$\Theta_l = d(\frac{\Theta_{l-1}}{s_0^{n_{l-1}}}) \cdot (s_0^{n_l -1}).$$
Then, by   Lemma \ref{jj} (3),
$\Theta_l$ is a meromorphic section of ${\cal O}_{\bar{P}_{k+l}V}(n_l)$
with at most log-poles along $D$
(here $s_0=s_0 \circ (\pi_{0,k+l-1})_*$).\\

By Lemma \ref{sat} we may identify $\bar{P}_kV$ as $P_kV_0 \times \bar{G}$.
Then we have:

$$
 \begin{array}{cccrl}
\co_{ \bpkl ~(V)} ~(-1) ~|_{ \bar X_{k + l} (f) }
& ~~~ \subset
~~ & \co_{ \bpkl ~(V) } ~(-1)
~~ & \stackrel{ \alpha }{ \lra }\ \ \ \ \
~~ &  \co_{ P_{k + l} ~V_{0} } ~(-1) \\
& & & & \\
\da
&
& \da
&
& \hspace{.8in}\da \\
& & & & \\
 \bar X_{k + l} ~(f)
& \subset
& \bpkl ~(V)
& =\ P_{k + l} ~V_{0} \times \ol{G} \ \ \
& \stackrel{ p_{1} }{ \lra }\ P_{k + l} ~V_{0}\ \  \\
& & & & \\
\da
&
& \da
& \da\ \ \
& \hspace{.8in}  \da \\
& & & \\
\bar X_{k} ~(f)
& \subset
& \ol{P}_{k} ~(V)
& =\ \ \ P_{k} ~V_{0} \times \ol{G} \ \ \
& \stackrel{ p_{1} }{ \lra }\ \ P_{k} ~V_{0} \ \\
& & & & \\
\da
&
& \da
& \da\ \ \
& \hspace{.8in} \da \\
& & & & \\
\hspace{-.6in} \Gamma \stackrel{f}{ \lra } \bar  X(f)
& \subset
& \ol{G}
& =\hspace{.8in} ~ \ol{G}\ \ \
& \hspace{.1in} \lra\ \ \ \ 0
 \end{array}
$$
where $0\in G \se D$ and $p_1$ is projection
to the first factor. Define,  for $l \geq 0$:
$$W_l = \{a \in G : (f+a)_{[k+l]}(u) \in \bar X_{k+l}(f) \; {\rm and} \;
\frac{\Theta_i}{s_0^{n_i}}|_{f_{[k+i]}(u)} =
\frac{\Theta_i}{s_0^{n_i}}|_{(f+a)_{[k+i]}(u)}, \; i=0,...,l \}.$$

\begin{lem} \label{dimV}
With the hypothesis and the setup as in the Main Lemma,
$W:= \bigcap_{l=0}^{\infty}W_l$ is an algebraic subvariety of $G$
and ${\rm dim}_0W \geq 1$.
\end{lem}

\noindent {\bf Proof of Lemma \ref{dimV} }
$W_l$ is an algebraic subvariety of $G$ as the group
action of $G$ on itself is algebraic.
Hence, $W$ is also algebraic.
Let $l_1 > l_2$ and $\pi_{k+l_2,k+l_1}: \bar{P}_{k+l_1}V \r \bar{P}_{k+l_2}
V$. If
$(f+a)_{[k+l_1]}(u) \in \bar X_{k+l_1}(f)$, then
$$(f+a)_{[k+l_2]}(u) = \pi_{k+l_2,k+l_1} \circ (f+a)_{[k+l_1]}(u)
\in \pi_{k+l_2,k+l_1}(\bar X_{k+l_1}(f)) = \bar X_{k+l_2}(f).$$
Hence, $W_l$, $l \geq 0$, is a decreasing sequence of algebraic subvarieties
of $G$. So the proof of Lemma \ref{dimV} is complete if we show:
$${\rm dim}_0W_l \geq 1 \; {\rm for} \; l \geq 0.$$

By the beginning of part iii) of the proof of Theorem 6.8 of Demailly in
\cite{Dem},
the rational map, obtained
by a basis of holomorphic sections of the line bundle ${\cal O}_{P_{k+l}V_0}
(2(k+l)-1)$, is a morphism on the subset $P_{k+l}(V)^{\rm reg}_0$ of regular
jets in $P_{k+l}V_0$ and seperates all points there.
Denote by $L_{(2(k+l)-1)}$ the linear system obtained by the pull backs of
these
sections by the map $\alpha$ (see the last diagram). Then we have
$$(L_{(2(k+l)-1)})^{\frac{n_l}{2(k+l)-1}} \s H^0(\bar{P}_{k+l}V,
{\cal O}_{\bar{P}_{k+l}V}(n_{l})).$$
Therefore,  the sections of $H^0(\bar{P}_{k+l}V, {\cal
O}_{\bar{P}_{k+l}V}(n_{l}))$
still separate points in the subset of regular jets of each
fiber of the map $\bar{P}_{k+l}V \r \bar{G}$.  Let the
map $\Phi_l: \bar{P}_{k+l}V \r \pr^{N_l}$ be defined by a basis of these
sections. Then the fiber of the map $\Phi_l$ through a regular jet $\xi \in
 \bar{P}_{k+l}V$ is necessarily
of the form $\{ \xi +a, \; a \in R \}$, where $R \s \bar{G}$ is algebraic.

To the basis of holomorphic sections
which define the map $\Phi_l$, we now add some extra
sections which we allow in addition to have log-poles along the divisor $D$,
namely the sections
$\Theta_i \cdot s_0^{n_l - n_i}$, $i=0,...,l$.  So we get a map
$$\tilde{\Phi_l}: \bar{P}_{k+l}V \r \pr^{N_l+l+1}.$$
This map will in general only separate  the subset of regular jets
of those fibers of the map $\pi_{0,k+l}: \bar{P}_{k+l}V \r \bar{G}$
 which are not over
the divisor $D$. But the fibers of the map $\tilde{\Phi}_l: \bar{P}_{k+l}V
\r \pr^{N_l+l+1}$ through a regular jet $\xi \in \bar{P}_{k+l}V \se
\pi_{0,k+l}^{-1}(D)$ must still be of the form $\{\xi +a, \; a \in R_\xi \}$,
where $R_\xi\s R \s \bar{G}$ is an algebraic subset.
So, Lemma \ref{dimV} is proved if we show that
$\tilde{\Phi_l}:\bar X_{k+l}(f) \r \pr^{N_l+l+1}$
has positive dimensional fiber through $\xi =f_{[k+l]}(u)$.

For proving this, we want to use our Ahlfors Lemma \ref{ahlfors}.
But this only applies for holomorphic sections.   We first extend the
diagram in equation (\ref{plp}) to
$$
 \begin{array}{rccccccl}
{ \bar Z}_{k + l}(f)
& \subset
& \bar P_{k + l} ~Y
& \stackrel{ p_{[k + l]} }{ \lra } \bar P_{k + l} ~G
& \supset
& \ol{P}_{k + l} ~(V)
& \supset
& X_{k + l}(f)  \\
& & & & & & & \\
&
& \da
& \da
&
& \da
& \\
& & & & & & & \\
\bar Z_{k}(f)
& \subset
& \bar P_{k} ~Y
& \stackrel{ p_{[k]} }{ \lra } \bar P_{k} ~G
& \supset
& \ol{P}_{k} ~(V)
& \supset
& X_{k}(f) \\
& & & & & & & \\
&
& \da
& \da
&
& \da
& \\
& & & & & & & \\
&
& (~ \ol{Y}, E ~)
& \stackrel{p}{\lra} (~ \ol{G}, S ~)
& \supset
& ( \bar G, S, \bar V )
&
 \end{array}
$$
where
 $${\bar Z}_k(f)=\overline{p_{[k]}^{-1}(X_k(f) \se \pi_{0,k}^{-1}
(D))}^{{\rm Zariski}} \s \bar Z_k \s \bar P_kY,$$
  $${\bar Z}_{k+l}(f)=\overline{p_{[k+l]}^{-1}(X_{k+l}(f)
 \se \pi_{0,k+l}^{-1}(D))}^{{\rm Zariski}} \s \bar Z_{k+l} \s \bar P_{k+l}Y.$$
By functoriality of the jet bundles
and definition  \ref{deflp}, the sections which define $\tilde{\Phi_l}$
pull back to holomorphic sections to span a linear system $$\tilde{L}_{k+l} \s
H^0\Big(\bar Z_{k+l}(f),{\cal O}_{\bar{P}_{k+l}Y}(n_{l})\Big).$$
The elements of
$\tilde{L}_{k+l}$ define the pull back of $\tilde{\Phi_l}$ to
$\bar Z_{k+l(f)}$ (which we denote again by $\tilde{\Phi_l}$).

So we can apply our Ahlfors Lemma \ref{ahlfors} to the map
$$\tilde{\Phi}_l : \bar Z_{k+l}(f) \r \pr^{N_l+l+1}$$
to conclude that $\~f_{[k+l]}(u) \in B_{k+l}(\bar Z_{k+l}(f))$, where
$\~f=p^{-1}\circ f$ and, without loss of generality,
$u\not\in {\rm Sing}\tilde{X}_{k+l}$.
For otherwise, the map $\~f$ would be constant or
(in the case of $\Gamma= \Delta^*$) at least extendable, which would
imply that $f$ has this property, or
the image of $\~f_{[k+l]}$
would be contained in ${\rm Sing}\bar Z_{k+l}(f)$, which is impossible,
since $\bar Z_{k+l}(f)$ is the proper transform of the Zariski closure of the
image of $f_{[k+l]}$.
Hence, $\tilde{\Phi}_l : \bar Z_{k+l}(f) \r \pr^{N_l+l+1}$
has positive dimensional fiber through $\~f_{[k+l]}(u)$.
Since $f_{[k+l]}(u) \not\in \pi_{0,k+l}^{-1}(S \cup D)$,
the map $p_{[k+l]}^{-1}$
is an isomorphism around $f_{[k+l]}(0)$. Hence,
$\tilde{\Phi_l}: X_{k+l}(f) \r \pr^{N_l+l+1}$ also has positive dimensional
fiber. This ends the proof of Lemma \ref{dimV}. \qed

\vspace{.3cm}

Let us now continue with the proof of the Main Lemma.
Without loss of generality we may assume that $f_{[k]}(u) \not\in
{\rm Sing}(X_k(f))$. Then there exists a neighborhood $U=U(0) \s G$
such that, for all $a \in U$, we have $(f+a)_{[k]}(u) \not\in
{\rm Sing}(X_k(f))$. By Lemma \ref{jj} (2), we get
$$((f+a)_{[k]})_{[l]}(u) \in \bar P_l(\tilde{X_k}(f))
\s \bar P_l({P}_k(V))$$ for $a \in W \cap U$ and all $l \geq 0$,
where we have omitted to write
the map $\Psi$. This is justified by the fact that around $(f+a)_{[k]}(u)$,
for $a \in U$, the variety $X_k(f)$ is smooth, and so $\Psi$ is an
isomorphism there.
Applying Lemma \ref{id}, we  get that $(f+a)_{[k]}(\Gamma) \s
\tilde{X_k}(f)$. Hence, $$(f+a)_{[k]}(\Gamma) \s X_k(f)$$ for $a \in W \cap U$.
But this means $$f(\Gamma) \s X_k(f) \cap (X_k(f) +a).$$ Since
$X_k(f)$ was the Zariski closure of $f(\Gamma)$, we get
$$X_k(f) = X_k(f)+a$$ for all $a \in W \cap U$.
We next want to show:
\begin{lem} \label{trans}
$$\left. \frac{\Theta}{s_0^{{\rm deg} \Theta}} \right|_{X_k(f)}$$
is invariant under the action of all $a \in W \cap U$
\end{lem}

\noindent {\bf Proof of Lemma \ref{trans} }
Let $a \in W \cap U$ be fixed. In order to simplify notation, we denote by
$F_{(i)}$, $i \in \nz_0$, the following rational function on $\bar{P}_{k+i}V$:
$$ F_{(i)} (y) := \frac{\Theta_i}{s_0^{n_i}} (y) -
\frac{\Theta_i}{s_0^{n_i}} (y+a).$$
Set $F=F_{(0)}$.
It suffices to show that, for all $i \in \nz_0$, we have
\begin{equation} \label{t0n}
\frac{\partial^i}{\partial t^i} (F \circ f_{[k]})(u) =0.
\end{equation}
For then, by analytic continuation applied to
$f_{[k]}: \Delta \r \bar{P}_kV$, we will
have $F \circ f_{[k]}(\Gamma) \equiv 0$,
so that $F \equiv 0$ on $X_k(f)$ as required
by Lemma \ref{trans}.

By abuse of notation, we have identified
$s_0 \circ (\pi_{0,k+l-1})_*$ with $s_0$.
Since these sections are maps from ${\cal O}_{\bar{P}_{k+l}V}(1)$ respectively
${\cal O}_{\bar{P}_{1}V}(1)$, this means we have
\begin{equation} \label{t1n}
s_0(f_{[k+l]}(t)) \cdot f_{[k+l-1]}'(t) = s_0(f_{[1]}(t)) \cdot f'(t)\: .
\end{equation}
Recall that the section $s_0$ is nonvanishing at $f_{[1]}(u)$, and that
$f'(u) \not= 0$. So, after possibly shrinking $\Delta$, we may reparametrize
$f$ such that
\begin{equation} \label{t2n}
s_0(f_{[k+l]}(t)) \cdot f_{[k+l-1]}'(t) = s_0(f_{[1]}(t)) \cdot f'(t)
\equiv 1\: .
\end{equation}
For the rest of this proof we fix the parameter $t$ in such a way that
equation (\ref{t2n}) is satisfied. We claim that
for any $i \in \nz_0$ we have:
\begin{equation} \label{t3n}
\frac{\partial^i}{\partial t^i} (F \circ f_{[k]})(t) = F_{(i)} \circ
f_{[k+i]}(t)\:.
\end{equation}

We prove this by induction over $i \in \nz_0$. The case {\bf S(0)} is clear
by definition. Assume that {\bf S($i$), $i <l$} is true. Then we have
\begin{eqnarray*} \frac{\partial^l}{\partial t^l} (F \circ f_{[k]})(t)
&=&  \frac{\partial}{\partial t}
    (\frac{\partial^{l-1}}{\partial t^{l-1}} (F \circ f_{[k]}))(t) \\
&=&  \frac{\partial}{\partial t}
     (F_{(l-1)} \circ f_{[k+l-1]})(t) \\
&=&  (dF_{(l-1)})((f_{[k+l-1]})(t)) \cdot  f_{[k+l-1]}'(t)\\
&=&  \frac{(dF_{(l-1)})((f_{[k+l-1]})(t)) \cdot
 f_{[k+l-1]}'(t)}{s_0((f_{[k+l]})(t)) \cdot  f_{[k+l-1]}'(t)}\\
&=&  \frac{(dF_{(l-1)})((f_{[k+l]})(t)) \cdot
 f_{[k+l-1]}'(t)}{s_0((f_{[k+l]})(t)) \cdot  f_{[k+l-1]}'(t)}\\
&=&  \frac{dF_{(l-1)}}{s_0}((f_{[k+l]})(t)) = F_{(l)}((f_{[k+l]}(t)).
\end{eqnarray*}
Here, we regard the differentials as linear maps on ${\cal O}(-1)$, more
specifically,
sections of ${\cal O}(1)$. To see the second equality from below,
recall that the section
$dF_{(l-1)}$ of ${\cal O}_{\bar P_1({P}_{k+l-1}(V))}(1)$
naturally restricts to a section of
${\cal O}_{\bar{P}_{k+l}V}(1)$.
This proves equation (\ref{t3n}). But from equation (\ref{t3n}),
equation (\ref{t0n}) follows immediately. This is because, by
the definition of $W$, $F_{(i)} \circ f_{[k+i]}(u) =0 \: $
for all $i \in \nz_0$. \qed

 Now the proof of the Main Lemma is immediate. Let us  define
$$\tilde{W} = \{ a \in G : X_k(f) = X_k(f)+a, \;
\frac{\Theta}{s_0^{{\rm deg}\Theta}} \; {\rm is\; invariant\; under}\; a \}.$$
$\tilde{W}$ is clearly a group, which is algebraic by Lemma \ref{sav}.
It is of positive dimension, since $(W \cap U) \s \tilde{W}$ and
${\rm dim}_0(W \cap U) \geq 1$. This finishes the proof of the Main Lemma
\ref{ji}. \qed

\subsection{Proof of Theorems \ref{c} and \ref{bp} }

{\bf Proof of  Theorems \ref{c} (a) and \ref{bp} (a) }
We apply the Main Lemma \ref{ji} in the special case where $k=0$ and
$V = \bar{T}G$. Then    Theorem \ref{bp} (a)
is immediate.  Theorem \ref{c} (a) is obtained by
dividing out by the biggest
algebraic subgroup of $G$ under which $X(f)$ is invariant. \qed $\ $

In order to prove the remaining parts of Theorems
\ref{c} and \ref{bp}, we first
have to choose the section $\Theta$ appropriately. We do this in the same way
as Siu-Yeung (\cite{SY}) or Noguchi (\cite{PNog}).

By Noguchi (\cite{PNog}, Lemma 2.1), there exists a theta function for
$D \s G$. This means the following.
Let $\pi : \cz^n \r G$ be the universal covering with a
`semi-lattice' $\Pi_1(G)$. There exists an entire function
$\theta : \cz^n \r \cz$ such that $$(\theta) = \pi^*D.$$
Moreover, for
any $ \gamma \in \Pi_1(G)$, there is an affine linear function
$L_{\gamma}$
in $x$ with
\begin{equation} \label{lgamma}
\theta (x+ \gamma) = e^{L_{\gamma}(x)} \theta (x), \; x \in \cz^n.
\end{equation}
But the proof which Noguchi gives actually yields more. Consider $G$ as a
$(\cz^*)^\ell$-principal fiber bundle over an abelian variety $A$,
and denote the projection map by $p: G \r A$.
Let $\pi : \cz^{\sf  m} \r A$ be the universal covering.
Then the fibered products of $\pi$ with $G$ respectively $\bar{G}$ over $A$ are
$(\cz^*)^\ell \times \cz^{\sf m}$ respectively $(\pr^1)^\ell \times
\cz^{\sf m}$. Let
$\tilde{\pi}: \cz^n \r (\cz^*)^\ell \times \cz^{\sf m}$ be the
universal covering. Then Noguchi's proof yields that there
exists a holomorphic function $\tilde{\theta}:
(\cz^*)^\ell \times \cz^{\sf m} \r \cz$ which extends to
a meromorphic function on $(\pr^1)^\ell \times \cz^{\sf m}$
such that $\theta = \tilde{\theta} \circ \tilde{\pi}$ is a
theta function for $D \s G$ as above. More precisely, if $(w_1,...,w_t)$ is
a multiplicative coordinate system of $(\cz^*)^\ell$
and $U \s \cz^{\sf m}$ is a small
neighborhood of a point in $\cz^{\sf m}$, then
$\tilde{\theta}|_{(\cz^*)^\ell \times U}$
 can be written as
\begin{equation} \label{thet}
 \sum_{{\rm finite}} a_{l_1...l_t}(y)w_1^{l_1}...w_t^{l_t},
\end{equation}
where the coefficients $a_{l_1...l_t}(y)$ are holomorphic functions on $U$.
As $(\pr^1)^\ell \times \cz^{\sf m}$ is the universal covering
space of $\bar G$, $\~\theta$ gives a multivalued defining
function for $D$ on $\bar G$ locally given by equation (\ref{thet}).
It follows that, whenever an algebraic expression in $\theta$
descends to $G$, it extends to $\bar{G}$ to a meromorphic object, and,
hence, to
a rational object.

Let $ \mu : \tilde{\Gamma} \r \Gamma$ be the universal cover, and let
$\tilde{f} : \tilde{\Gamma} \r \cz^n$ be any lift of the map
$f \circ \mu : \tilde{\Gamma} \r G \se D$.  Then we can choose a linear
coordinate
system $(z_1,...,z_n)$ of $\cz^n$ such that in these coordinates, $\tilde{f}$
is
 expressed as
$$\tilde{f}(z) = (f_1(z),...,f_{n'}(z),0,...,0)$$ with holomorphic functions
$f_1 (z),...,f_{n'}(z)$, for which the functions\\ $1,f_1,f_2,...,f_{n'}$
are linearly independent.

The set of differential equations $dz_{n'+1}=0$,...,$dz_n=0$ obviously
defines a subbundle $\tilde{V} \s T\cz^n$,
 which is invariant under translation and, hence,
descends to  $G$.
 It is important to remark that this subbundle extends
to a special subbundle $V \s \bar{T} G$. This is true,
 since by definition, $V$ is
just of the form $V=\cz^{n'} \times G$ with respect to the
trivialization of $\bar{T} G$ given by the standard logarithmic forms
$dz_1,...,dz_n$. \\

Siu-Yeung (\cite{SY}) defined the following logarithmic jet differential:
$$ \Theta = \left|
\begin{array}{llll}
d \log \theta & dz_1 & ... & dz_{n'}\\
d^2 \log \theta & d^2 z_1 & ...& d^2 z_{n'}\\
. & . & . & . \\
. & . & . & . \\
. & . & . & . \\
d^{n'+1} \log \theta & d^{n'+1}z_1 & ... & d^{n'+1}z_{n'}
\end{array}
\right| $$
We want to use this $\Theta$ as the
$\Theta$ in our Main Lemma \ref{ji}. So we need:

\begin{lem} \label{theta}
$\Theta$ can be considered as a meromorphic section in the line
bundle $${\cal O}_{\bar{P}_{n'+1}V}(\frac{(n'+1)(n'+2)}{2})$$
with at most log-poles along $D$.
\end{lem}

\noindent {\bf Proof of Lemma \ref{theta} } We first show (part (a))
that $\Theta$ is a meromorphic function on
$J_k\cz^n$, and, hence, also on $J_k\~V$,
which satisfies equation (\ref{invj}) with order $k=n'+1$
and degree $m={(n'+1)(n'+2)}/{2}=k(k+1)/2$
(recall that $\tilde{V} \s T\cz^n$ is defined by
$dz_{n'+1}=...=dz_n=0$).
By using equation (\ref{thet}) it follows actually that $\Theta$ defines such a
function on $J_k (\cz^{*\ell}\times \cz^{\sf m})$,
which extends, again by equation (\ref{thet}), to a meromorphic function on
$\bar J_k (\cz^{*\ell}\times \cz^{\sf m})$. This meromorphic
function descends to a multivalued function on
$\bar J_k G$.
Then we show (part (b))  that it is actually singlevalued
on $\bar{J}_k V \s \bar J_k G$.
So it  corresponds, by Proposition \ref{dlog}, to a meromorphic section of
 ${\cal O}_{\bar{P}_{k}V}(m)$.
 That it is a meromorphic section of
 ${\cal O}_{\bar{P}_{k}V}(m)$ with
at most log-poles along $D$ then follows
from the local description of $\tilde{\theta}$
in equation (\ref{thet}) and  from the local
description of meromorphic sections with
at most log-poles along $D$ right after definition
\ref{deflp}: In fact, by this description, applied to the multivalued
meromorphic fuction on $\bar J_k G$, it gives rise to a (possibly) multivalued
holomorphic function on $\bar J_k Y$. But, as we saw above, it
is singlevalued over $\bar{J}_k V$, and hence, it is also
singlevalued over $\bar J_k Y$. Thus, it yields a meromorphic
section with at most log-poles along $D$. \\

(a) We show more generally: Let $h_1,...,h_r: (\Delta, 0) \r \cz$ be
nonvanishing germs of holomorphic
functions.
Let $g_i = \log h_i$, $i=1,...,r$.
Then
$$\left|
\begin{array}{llll}
d g_1 & d g_2 & ... & d g_r\\
d^2 g_1 & d^2 g_2 & ... & d^2 g_r\\
. & . & . & . \\
. & . & . & . \\
. & . & . & . \\
d^{r} g_1 & d^{r} g_2 & ... & d^{r} g_r
\end{array}
\right| $$
gives a jet differential on $J_r \Delta_{0}$
which is equivariant under the full reparame\-tri\-zation group $J_r \Delta_0$
in the sense of equation (\ref{invj}).
(We will apply this for the case where the germs $h_j$ are obtained by
composing $\theta$
and the exponentials of the $z_i$'s with the germ
$f\colon (\cz,0)\r \cz^n$ representing
the jet in our case.)
\\

 Only the equivariance is nontrivial. Let $g=(g_i)$.
Let $\phi \in J_r\cz_0$.
Then, by using the identity
$$ (g \circ \phi)^{(j)}(0) = g^{(j)}(0) (\phi'(0))^j +
\sum_{s=1}^{j-1} \sum_{i_1+...+i_s=j} c_{i_1...i_s} g^{(s)}(0)
\phi^{(i_1)}(0)...\phi^{(i_s)}(0),$$
we get by induction on $r$ that
$$(g \circ \phi)' \wedge ... \wedge (g \circ \phi)^{(r)}(0)
= g' \wedge ...\wedge g^{(r)} \cdot (\phi'(0))^{\frac{r(r+1)}{2}}.$$
This gives the desired equivariance.\\

(b) {From} equation (\ref{lgamma}) we have, for $\gamma \in \Pi_1(G)\s \cz^n$:
$$d^i \log \theta (x+ \gamma) = d^i \log \theta (x) + d^i L_{\gamma}(x)
= d^i \log \theta (x) + \sum_{j=1}^{n'} a_jd^ix_j +
\sum_{j=n'+1}^n a_j d^ix_j,$$
where $a_i\in \cz$ are constants.
Then, from the properties of the determinant and the fact that we restrict
$\Theta$ to $J_{n'+1} \tilde{V}$, it follows that this jet differential is
invariant under the action of $\Pi_1(G)$.
Hence, it descends to $G$. \qed

\begin{lem} \label{S}
Under the assumptions of Theorem \ref{c}
  (b), or Theorem \ref{bp} (b) and (c) and the additional assumption that
$f$ does
not extend, the following holds:\\
 If $X(f)\cap D \not= \emptyset$, then $X(f) \cap D$ is foliated by translates
of an algebraic subgroup $G'' \s G'$ of positive dimension, where $G'$
is the maximal subgroup whose translates foliate $X(f)$.
\end{lem}

\noindent {\bf Proof of Lemma \ref{S} }
We may assume that $f$ is nonconstant and,
for the case that $\Gamma = \Delta^*$, that  the map $f$ does not extend.
We apply the Main Lemma \ref{ji} to get the existence of an algebraic subgroup
$G'' \s G$ of positive dimension which leaves $X_{n'+1}(f)$ and
$\Theta|_{X_{n'+1}(f)}$
invariant. As $X_{n'+1}(f)$ is invariant under the action
of $G''$ and the projection
$\pi_{1,n'+1}$ (respectively $\pi_{0,n'+1}$) maps $X_{n'+1}(f)$ surjectively
onto $X_1(f)$ (respectively $X(f)$),
we see that $X_1(f)$ and $X(f)$ are also invariant.

Take any $a \in G''$. Since $\Theta |_{X_{n'+1}(f)}$ is invariant under
translation by $a$,
$$\left|
\begin{array}{llll}
d \log \frac{\theta (f)}{\theta (f+a)} & df_1 & ... & df_{n'}\\
d^2 \log \frac{\theta (f)}{\theta (f+a)} & d^2 f_1 & ... & d^2 f_{n'}\\
. & . & . & . \\
. & . & . & . \\
. & . & . & . \\
d^{n'+1} \log \frac{\theta (f)}{\theta (f+a)} & d^{n'+1}f_1 & ...
& d^{n'+1}f_{n'}
\end{array}
\right| \equiv 0 \; {\rm on} \; \Gamma .$$
Now $d \log \frac{\theta (x)}{\theta (x + a)}$
is a rational differential on $G$. Since $f+a$ cannot map entirely into the
zero set of
$\theta$, because $X(f)$ is the Zariski closure of $f(\Gamma)$,
$$\frac{\partial}{\partial z}(\log \frac{\theta (f)(z)}{\theta (f+a)(z)}),
\frac{\partial}{\partial z}f_1(z),...,\frac{\partial}{\partial z}f_{n'}(z)$$
are well defined meromorphic functions on $\Gamma$.
The functions $\frac{\partial}{\partial z}f_1(z),...,$
$\frac{\partial}{\partial z}f_{n'}(z)$ are linearly
independent as $1,f_1,...,f_{n'}$ were so. Hence, we get,
 by the classical Lemma of the Wronskian \cite{Cartan}, that there exist
 complex numbers
$c_1,...,c_{n'}$ (which may depend on $a \in G''$) such that
$$ d  \log \frac{\theta (f)(z)}{\theta (f+a)(z)} + c_1df_1(z) + ... +
c_{n'}df_{n'}(z) \equiv 0 \; {\rm on} \; \Gamma.$$
So we have
\begin{equation} \label {*}
d \log \frac{\theta (x)}{\theta (x+a)} + c_1dx_1 + ... +
c_{n'}dx_{n'} \equiv 0
\end{equation}
on $f_{[1]}(\Gamma)$. Moreover, since $d \log \frac{\theta (x)}{\theta (x +
a)}$
is a rational differential on $G$, this equation
holds on $X_1(f)$.\\

Assume now that Lemma \ref{S} does not hold. Then
there exists $x_0 \in X(f) \cap D$ and $a_0 \in G''$
such that $x_0+a_0 \not\in D$. We want to show that this assumption leads
to a contradiction.  {From} equation (\ref{*}) we get that
\begin{equation} \label{**}
d \log \frac{\theta (x)}{\theta (x+a_0)}
= d \log  \frac{\theta (x+b)}{\theta (x+a_0+b)}
\end{equation}
on $X_1(f)$ for $b \in G''$. This means that
$$d \log (\frac{\theta (f)}{\theta (f+b)}
  \frac{\theta (f+a_0+b)}{\theta (f+a_0)}) \equiv 0 \; {\rm on} \; \Gamma .$$
Hence,
$$ \frac{\theta (f)}{\theta (f+b)}
  \frac{\theta (f+a_0+b)}{\theta (f+a_0)} = c_{a_0,b}\; {\rm on} \; \Gamma, $$
where $c_{a_0,b} \in \cz$ is a constant, which may depend on $a_0$ and $b$.
Since $\frac{\theta (x)}{\theta (x+b)}
  \frac{\theta (x+a_0+b)}{\theta (x+a_0)}$
  is a well defined rational function on $G$, we have
  \begin{equation} \label{***}
\frac{\theta (x)}{\theta (x+b)}
  \frac{\theta (x+a_0+b)}{\theta (x+a_0)} = c_{a_0,b} \;  \:{\rm on}\ X(f),
  \end{equation}
where $b \in G''$.
Now $x_0+a_0 \not\in D$, but $x_0 \in D$. So we get,
for $b = a_0$ and $x=x_0$, that
$c_{a_0,a_0}=0$. This means, as $X(f)$ is irreducible, that either
$\theta (x) \equiv 0$ or $\theta (x +2a_0) \equiv 0$ on $X(f)$. But both
is not true, as one sees by taking $x=x_0+a_0$ respectively
$x=x_0+a_0-2a_0$ (remark
that the latter is still in $X(f)$, since $X(f)$ is invariant under the action
of $G''$). So our assumption was wrong, and we have proved Lemma \ref{S}. \qed
\vspace{.3cm}

\noindent {\bf Proof of Theorem \ref{c} (b) }
Assume that $X(f) \cap D \not= \emptyset$. We want to show that this assumption
leads to a contradiction. After applying a translation, we may assume, by
Theorem \ref{c} (a), that $X(f)$ again is a semi abelian variety $G'$ with
nonempty divisor
$D'$ in $G'$, where $D'$ is the reduction of $X(f) \cap D$.  Now we devide
through
the maximal algebraic subgroup $\tilde{G}$ of $G'$ which foliates $D'$. Then,
 by applying
Lemma \ref{S} to the quotient $G' / \tilde{G}$ and by taking the invers
image under the quotient map,
we get $X(f) \cap D = \emptyset$, which contradicts our assumption. \qed
\vspace{.3cm}

\noindent {\bf Proof of Theorem \ref{bp} (b) and (c) }
 Part (b) follows immediately from Lemma  \ref{S}.
For (c), let $G$ be an abelian variety.
Let $G'$
again be the maximal algebraic subgroup the translates of which foliate $X(f)$.
We may assume that all translates
of $G'$ which foliate $X(f)$ intersect $D$ (in particular $X(f) \cap D
\not= \emptyset$), for otherwise we finish the
proof by using Lemma \ref{proj}. Then there must be such
a translate $T_0$ of $G'$ such that $T_0 \not\s D$. Now by Lemma \ref{S}
we find a subgroup $G'' \s G'$ of positive dimension which foliates
$X(f) \cap D$. Hence, $T \cap D$ is foliated by translates
of $G''$. But since $T$ is also foliated by translates of $G''$, there must
be such a translate not hitting $D$ at all. This finishes the proof again
by using Lemma \ref{proj}. \qed

\section{Appendix}
\setcounter{equation}{0}

We use the notations of subsection 4.1.
 We now give a key result via which pseudometrics of negative curvature
are usually constructed (see (2.7) of \cite{GG}).
We point out that our version is sharper than the ones in \cite{Dem, GG} for
the basic locus in our definition is smaller than theirs.

\begin{lem}\label{Bs}
 Let the setup be as in subsection 4.1, and assume further that $X$ is
normal.  Then given
line bundles $L$ and $H$ over $X$, there is an integer
$l _{1} \geq 1$  such that $x \notin E_{L}$ implies that
$x \notin {\rm Bs} |l L -H |$  for all positive multiples $l$
of $l _{1}$, more specifically,
$$
E _L \supseteq {\rm Bs} | l L - H | .
$$
\end{lem}

\noindent {\bf Proof }  Observe that we may always write
$H + H ^{\prime} = H ^ {\prime \prime}$, where
$H ^{\prime}$ and $H ^{\prime \prime}$ are very ample divisors.  Then
$\mbox{Bs} | l L - H ^{\prime \prime} | \supseteq \mbox{Bs} | l L - H |$,
as one can see from the fact that
$
\mbox{Bs} |G-G ^{\prime} | \cup \mbox{Bs} | G ^{\prime} | \supseteq \mbox{Bs}
|G|
$
for arbitrary line bundles $G$ and $G^{\prime}$.
Hence, we will assume without loss of generality that $H$ is very ample.

Let $x \in X$ be outside $E _{L}$.
Then,  we may assume
that $\ph _{L}$ is birational onto
its image after replacing $L$ by a suitable multiple of $L$
(see 1.10 and 5.7 of \cite{Ueno}).
It will be
sufficient to show that $x$ is outside ${\rm Bs}|lL-H|$ for some $l$, and
hence, for all multiples thereof, as  Lemma \ref{Bs} would then
follow by the quasi-compactness of $X \setminus E _{L}$.

Consider the ideal sheaf ${\cal I} \supseteq {\cal O} _{X}$ generated by
the global sections of $L$.  By blowing up this ideal sheaf,
we obtain a
modification $\sigma : \~X \r X$ so that ${\cal J} = \sigma ^{\star}
{\cal I} \cdot {\cal O} _{\~X}$ is an invertible ideal sheaf of ${\cal
O}_{\~X}$
generated by a global section $s$ of the line bundle
$F = {\cal O}_{\sigma} (-1)$, namely,
$
{\cal J} = \mbox{Im} \{ {\cal O} (F ^* )
\stackrel{\otimes s}{\longrightarrow} {\cal O} \} .
$
Then
$
\bar{L} := \sigma ^ {-1} L  - F
$
is spanned by the sections $(\sigma ^ {-1} t)/s$ as $t$
ranges over $H^{0} (L)$.  Note that $E_{\bar L}=\sigma^{-1}(E_L)$,
that $\sigma^{-1}$ is an isomorphism on the neighborhood $X\setminus E_L$
of $x$ and that any section of
$l\bar L-\sigma^{-1}H$
not vanishing on the point $\sigma^{-1} (x)$
gives rise to a section of $lL-H$ not vanishing on $x$ by tensoring
with $s^l$ (Zariski's Main Theorem and $s(x)\ne 0$). Hence, replacing $(X,L)$
by $(\~X,\bar L)$ we may assume that $\ph _{{L}}:X \r \pr ^{n}$
is a birational morphism onto its image $W=\ph_L(X)$.
Let
$
\sigma _{0} : W _{0} \r W \subseteq \pr ^{n}
$
be the normalization of $W$. Then $H _{0} := \sigma _{0} ^ {-1}
{\cal O} _{\pr ^{n}} (1)$ is ample so that
there is a positive integer $d$ such that $d H_{0}$ is very ample on
$W _{0}$.  As $X$ is normal, there is a canonical morphism
$\ph : X \r W_{0}$
such that $\ph _{{L}} = \sigma _{0} \circ \ph$. Noting
$
\ph ^ {-1} H _{0} = {L},
$
we see that the image $W _{1}$ of the morphism $\ph _{d {L}}$ admits a
birational morphism $r$ to $W _{0}$ and that $\ph=r\circ\ph_{dL}$.
As $\ph$ is connnected by Zariski's Main Theorem, $\ph^{-1}(\ph (x))=\{x\}$.
Hence, replacing ${L}$ by $d {L}$ we may
assume that
\begin{equation}\label{1-1}
\ph _{{L}} ^ {-1} (\ph _{{L}} ( x )) = \{ x \}.
\end{equation}

As $H$ is very ample, $|H|$ has an element
$D$ such that $x \notin
D \not\subseteq E _{{L}}$
by general positioning.
We may now choose, thanks to equation (\ref{1-1}),
a hypersurface of sufficiently high degree
$l$ in $\pr ^{n}$ containing $\ph _{{L}} (D)$ but not
$\ph _{{L}} (x )$.
This gives a divisor in $| l {L} - H |$ not containing $x$
as desired.  \qed

In practice, Lemma \ref{Bs}  is all that one uses.
But one can easily deduce the following strengthened version in order to
complete the picture.  \begin{lem}\label{S_L} Let $X$ be a normal complex
projective variety with
any line bundle $H$.  For any line bundle $L$  over
$X$, there is
an integer $m_0$  such that
$$
E _{m _{0} L} \supseteq S _{L} \supseteq \mbox{Bs} | m _{0} L - H | \, .$$
Moreover, if $H$ is very ample, the both inclusions are equalities.
\end{lem}

\noindent {\bf Proof }
Clearly there is an integer $N$ such that
$
S _{L} = \bigcap _{m=1} ^ {N} E _{mL} .
$
For each $m$, there is an integer $l _{m} > 0$, such that $E _{mL} \supseteq
\mbox{Bs} | l L - H|$ for all positive multiple $l$ of $l _{m}$ by  Lemma
\ref{Bs}.
 Letting $m _{0}$ be a common multiple of $l _{1} ,..., l _{N}$,
we see that
$
E _{m _{0} L} \supseteq S _{L} \supseteq \mbox{Bs} | m _{0} L - H |.
$
If, furthermore, $H$ is very ample, one easily verifies that $\mbox{Bs} | m
L - H | \supseteq E _{m L}$
for all $m$. \qed

\noindent{\bf Remark } As $\mbox{Bs}|lL-H|\supseteq E_{lL}\supseteq S_L$, it
follows that $S _{L} = \bigcap _{l > 0} \mbox{Bs} | lL - H|$ for
any very ample $H$.

{\small 

\noindent Gerd Dethloff\\
D\'{e}partement de Math\'{e}matiques,
UFR Sciences et Techniques\\
Universit\'{e} de Bretagne Occidentale\\
 6, avenue Le Gorgeu, BP 452\\
 29275 Brest Cedex, {France}.\\
e-mail: dethloff@univ-brest.fr\\

\noindent Steven Lu\\
Department of Pure Mathematics\\
 University of Waterloo\\
Waterloo, Ontario,
Canada N2L 3G1\\
e-mail: slu@math.uwaterloo.ca}

\end{document}